\RequirePackage{fix-cm}
\RequirePackage{amsmath}
\RequirePackage[hyphens]{url}
\documentclass[smallcondensed]{svjour3hack}       
\smartqed  

\usepackage{graphicx}
\usepackage{amssymb}
\usepackage{placeins}
\usepackage{xcolor}
\usepackage{subfig}
\usepackage{multirow}
\usepackage[page]{appendix} 

\usepackage{tikz}
\usetikzlibrary{fit}
\usetikzlibrary {arrows.meta}
\usetikzlibrary{shapes.geometric}  

\usepackage[linesnumbered, boxed, english, onelanguage, vlined]{algorithm2e}
\usepackage{mathtools}
\usepackage{mathrsfs}

\usepackage{amsmath, amsfonts, amssymb}
\usepackage{pdfpages}

\usepackage[normalem]{ulem} 

\usepackage[hidelinks,hypertexnames=false,colorlinks=true,breaklinks=true,bookmarks=true,urlcolor=blue,citecolor=blue,linkcolor=blue,bookmarksopen=false,draft=false]{hyperref}

\makeatletter
\let\c@proposition\c@theorem
\let\c@corollary\c@theorem
\let\c@lemma\c@theorem
\let\c@definition\c@theorem
\let\c@example\c@theorem
\let\c@remark\c@theorem
\let\c@jremark\c@theorem
\let\c@jconjecture\c@theorem
\makeatother

\makeatletter
\renewcommand*{\thetable}{\arabic{table}}
\renewcommand*{\thefigure}{\arabic{figure}}
\let\c@table\c@figure
\makeatother

\allowdisplaybreaks

\vbadness=\maxdimen

\makeatletter
\renewcommand*{\top}{%
  {\mathpalette\@transpose{}}%
}
\newcommand*{\@transpose}[2]{%
  \scriptsize
  \raisebox{\depth}{$\m@th#1\mathsf{T}$}%
}
\makeatother

\DeclareMathOperator{\Diag}{Diag}
\DeclareMathOperator{\diag}{diag}
\DeclareMathOperator{\ldet}{ldet}

\DeclareMathOperator{\Tr}{Tr}
\DeclareMathOperator{\rank}{rank}

\DeclareMathOperator{\argmin}{argmin}

\renewcommand{\circeq}{\mathrel{\ooalign{\hss$\circ$\hss\cr$\equiv$}}}

\newcommand{\znaturalthing}{\mathfrak{z}_{\mbox{\protect\tiny $\mathcal{N}$}}}
\newcommand{\znatural}{\hyperlink{znaturaltarget}{\znaturalthing}}

\newcommand{\zrelax}{\mathfrak{z}_{\mbox{\protect\tiny $\mathcal{R}$}}}

\newcommand{\zdoptthing}{{\mathfrak{z}}_{\text{\normalfont\protect\tiny  D-Opt}}}
\newcommand{\zdopt}{\hyperlink{zdopttarget}{\zdoptthing}}

\newcommand{\zmespthing}{z_{\mbox{\normalfont\protect\tiny MESP}}}

\newcommand{\znlpthingmesp}{{z}_{\text{\normalfont\protect\tiny NLP}}}
\newcommand{\znlpmesp}{\hyperlink{znlptargetmesp}{\znlpthingmesp}}

\newcommand{\znlpidthingmesp}{z_{\mbox{\normalfont\protect\tiny NLP-Id}}}
\newcommand{\znlpidmesp}{\hyperlink{znlpidtargetmesp}{\znlpidthingmesp}}

\newcommand{\znlpdithingmesp}{{z}_{\mbox{\normalfont\protect\tiny NLP-Di}}}

\newcommand{\znlptrthingmesp}{{z}_{\mbox{\normalfont\protect\tiny NLP-Tr}}}

\newcommand{\zlinxthing}{z_{\mbox{\protect\tiny linx}}}
\newcommand{\zlinx}{\hyperlink{zlinxtarget}{\zlinxthing}}

\newcommand{\zspectralthing}{z_{\mbox{\protect\tiny $\mathcal{S}$}}}
\newcommand{\zspectral}{\hyperlink{zspectraltarget}{\zspectralthing}}

\newcommand{\zhadamardthing}{z_{\mbox{\protect\tiny $\mathcal{H}$}}}
\newcommand{\zhadamard}{\hyperlink{zhadamardtarget}{\zhadamardthing}}

\newcommand{\zdiagonalthing}{z_{\mbox{\normalfont\protect\tiny diag}}}
\newcommand{\zdiagonal}{\hyperlink{zdiagonaltarget}{\zdiagonalthing}}

\newcommand{\zspectraldoptthing}{{\mathfrak{z}}_{\mbox{\protect\tiny $\mathcal{S}$}}}
\newcommand{\zspectraldopt}{\hyperlink{zspectraldopttarget}{\zspectraldoptthing}}

\newcommand{\zhadamarddoptthing}{{\mathfrak{z}}_{\mbox{\protect\tiny $\mathcal{H}$}}}
\newcommand{\zhadamarddopt}{\hyperlink{zhadamarddopttarget}{\zhadamarddoptthing}}

\newcommand{\zgammathing}{z_{\mbox{\protect\tiny $\Gamma$}}}

\newcommand{\zgammaplusthing}{z_{\mbox{\protect\tiny $\Gamma^{+}$}}}

\newcommand{\zbqpthing}{z_{\mbox{\normalfont\protect\tiny BQP}}}

\newcommand{\DOPT}{
\mbox{\rm D-Opt}}
\newcommand{\MESP}{
\mbox{\rm MESP}}

\newcommand{\NLPIdsym}{
\mbox{\rm{NLP-Id}}}
\newcommand{\NLPDisym}{
\mbox{\rm NLP-Di}}
\newcommand{\NLPTrsym}{
\mbox{\rm NLP-Tr}}

\newcommand{\Identity}{\hyperlink{Identitytarget}{\mbox{Ident}}}
\newcommand{\Diagonal}{\hyperlink{Identitytarget}{\mbox{Diag}}}
\newcommand{\Trace}{\hyperlink{Identitytarget}{\mbox{Trace}}}

\newcommand{\NLPId}{\hyperlink{NLPIdtarget}{\NLPIdsym}}
\newcommand{\NLPDi}{\hyperlink{NLPDitarget}{\NLPDisym}}
\newcommand{\NLPTr}{\hyperlink{NLPTrtarget}{\NLPTrsym}}

\begin{document}

\title{On the relationship between MESP and 0/1 D-Opt and their upper bounds\thanks{M. Fampa was supported in part by CNPq grant 307167/2022-4.
J. Lee was supported in part by AFOSR grant FA9550-22-1-0172.}}
\titlerunning{On MESP and 0/1 D-Opt}        
\author{Gabriel Ponte \and Marcia Fampa \and Jon Lee}%
\authorrunning{Ponte, Fampa \& Lee} 

\institute{           G. Ponte and J. Lee  \at
            University of Michigan.
              \email{gabponte@umich.edu}, \email{jonxlee@umich.edu}
\and
M. Fampa\at
            Federal University of Rio de Janeiro.
               \email{fampa@cos.ufrj.br}
}

\date{}

\maketitle
\begin{abstract}
We establish strong connections between two fundamental nonlinear 0/1 optimization problems coming from the area of experimental design, namely maximum-entropy sampling and 0/1 D-optimality. The connections are based on maps between instances,
and we analyze the behavior of these maps. 
Using these maps, we transport basic upper-bounding methods between these two problems.
Further, we establish results relating how different branch-and-bound schemes based on these maps compare.
Additionally, we observe some surprising numerical results, where bounding methods that did not seem promising in their direct application to real-data MESP instances, are now useful for MESP instances that come from 0/1 D-optimality.
\end{abstract}

\widowpenalty=10000
\clubpenalty=10000


\section{Introduction}\label{sec:int}
Some nice families of integer nonlinear-optimization problems come from the area of experimental design. One important 
problem is the Gaussian case of the maxi\-mum-entropy sampling problem (MESP). Here, we have an input covariance matrix of order $n$, and we wish to select a principal submatrix of order $s$, so as to maximize the 
``differential entropy'' (see, e.g., \cite{FLbook}).
In the Gaussian case, the differential entropy is, up to constants, the logarithm of the determinant.
 Another fundamental problem is the Gaussian case of the
0/1 D-optimality problem (0/1 $\DOPT$).
Briefly, the problem aims to select a subset of 
$s$ design points, from a universe of
$n$ given design points in $\mathbb{R}^m$, with the goal of 
minimizing the ``generalized variance'' of the least-squares 
parameter estimates; see, e.g., \cite{PonteFampaLeeMPB}
and the references therein. 
In our setup, we assume that we incorporate a further given set of $q$ design points in $\mathbb{R}^m$.
A notable special case  is the 0/1 D-optimal data-fusion problem, which assumes 
that the further $q$ design points span $\mathbb{R}^m$
(see, e.g., \cite{li2022d}).

\medskip
\noindent{\bf Brief literature review.}   $\MESP$  is a very-well studied topic  
 in experiment\-al design. A recent comprehensive reference on the state-of-the-art for branch-and-bound (B\&B) approaches is \cite{FLbook}, with references therein to 
background and previous work. Basic upper bounds are the ``spectral bound''
(see \cite{KLQ}), the ``diagonal bound'' (see \cite{HLW}), and upper bounds based on convex relaxation: the  ``factorization bound'' (see 
\cite{Nikolov,Weijun,FLbook,FactPaper}), the ``linx bound'' (see \cite{Kurt_linx}), the ``BQP bound'' (see \cite{Anstreicher_BQP_entropy}) and the ``NLP bound'' (see \cite{AFLW_Using}).
See \cite{FL_update} for some of the latest advances.

Integer D-Optimality, is a closely related 
and fundamental problem in experimental design. There are several variations, but we concentrate on the 0/1 version. A recent reference on the state-of-the-art for B\&B approaches is \cite{PonteFampaLeeMPB}, with  references therein to 
background and previous work. \cite{Welch} was the first to approach solving integer $\DOPT$ with an exact B\&B algorithm,
employing an upper bound based on Hadamard's inequality and another based on continuous relaxation (closely related to what we call these days the ``natural bound'' (see \cite{PonteFampaLeeMPB}, e.g.)). 
\cite{KoLeeWayne} proposed a spectral bound
and analytically compared it with a Hadamard bound
for 
0/1 D-optimal data-fusion instances only. \cite{li2022d} applied a local-search procedure and an exact B\&B algorithm to the 
0/1 D-optimal data-fusion problem. 
\cite{PFLXadmm} gave fast ADMM algorithms for convex relaxation of 0/1 $\DOPT$ (and $\MESP$). 

\smallskip
\noindent{\bf Organization and contributions.} 
In \S\ref{sec:equivalence}, we analyze four useful mappings, $\mathcal{M,P,D,F}$, between instances of $\MESP$ and instances of 0/1 $\DOPT$. \cite{li2022d} had worked extensively with $\mathcal{P}$ and $\mathcal{F}$ previously, establishing a close connection between the 0/1 D-optimal data-fusion problem (a special case of 0/1 $\DOPT$) and ``positive definite MESP''. 
We draw attention to the map $\mathcal{M}$, discussed mainly for special cases in \cite{PonteFampaLeeMPB}, which applies to arbitrary instances of 0/1 D-Opt. 
Finally, we introduce a completely new map $\mathcal{D}$, which applies to arbitrary
instances of $\MESP$.
This map allows us to re-cast, with computational advantage, the so-called ``NLP-Id'' upper bound for an arbitrary MESP instance as a ``natural bound'' calculation for an associated 0/1 D-Opt instance; see \S\ref{sec:compare-bounds}.
Further in \S\ref{sec:compare-bounds}, based on the equivalences established in \S\ref{sec:equivalence}, we can transfer any basic upper-bounding method for $\MESP$ instances to 0/1 $\DOPT$ instances, and conversely.
At present, there are six known basic upper-bounding methods for $\MESP$ and three for 
0/1 $\DOPT$.
In \S\ref{sec:bb}, in the context of the maps that we present between $\MESP$ instances and 
0/1 $\DOPT$ instances, we analyze the 
strength of upper bounds for child subproblems, and we establish 
the surprising relative strength of a particular B\&B strategy. 
In \S\ref{sec:numexp}, we present results of numerical experiments, demonstrating the benefits of applying various upper bounds.  
In particular, we explore the transference of upper bounds between $\MESP$ and 0/1 D-Opt, investigating the bounds transferred from $\MESP$ to 0/1 D-Opt via $\mathcal{M}$, the so-called ``$\mathcal{M}$-induced upper bounds for 0/1 D-Opt''; and from 0/1 D-Opt to $\MESP$ via $\mathcal{D}$, the so-called ``$\mathcal{D}$-induced upper bounds for $\MESP$''. We identify gains in both bound quality and computation time with the transference of bounds. Considering $\MESP$ instances derived by applying  $\mathcal{M}$ to instances of 0/1 D-Opt, we observe unusual comparative performance of the $\MESP$ upper bounds relative to those in the literature, interestingly demonstrating that new covariance matrix structures can alter the choice of the best bounds.

\smallskip
\noindent{\bf Notation.}
Throughout, we denote any all-zero square matrix by $0$, and we denote any all-zero (column) vector by $\mathbf{0}$.
We denote any all-one vector
by $\mathbf{e}$, any \hbox{$i$-th} standard unit vector by $\mathbf{e}_i$\,,  
 and the identity matrix of order $n$ by $I_n$\,.
 We let $\mathbb{S}^n$  (resp., $\mathbb{S}^n_+$~, $\mathbb{S}^n_{++}$)
 denote the set of symmetric (resp., positive-semidefinite, positive-definite) matrices of order $n$.
 We let $\Diag(x)$ denote the $n\times n$ diagonal matrix with diagonal elements given by the components of $x\in \mathbb{R}^n$, and we let $\diag(X)$ denote the $n$-vector with elements given by the diagonal elements of $X\in\mathbb{R}^{n\times n}$. 
  For matrix $X$,  $X_{S,T}$\, is the submatrix  with row (column) indices $S$ ($T$);
 if $S$ ($T$) is $N$, we write $X_{\cdot T}$ ($X_{S\cdot}$); we write $i$ for $\{i\}$.
We let
$\ldet$ denote the natural logarithm of the determinant. 
We let $\Tr$ denote the trace. 
For matrices of the same shape,
$X_1\circ X_2$ is the Hadamard (i.e., element-wise) product.
We let $N:=\{1,2,\ldots,n\}$.  
For $X\in\mathbb{S}^n$ and $i\in N$, we let $\lambda_i(X)$ denote its $i$-th greatest eigenvalue of $X$, and $\delta_i(X)$ denote its $i$-th greatest diagonal component. We let  $\lambda(X):=(\lambda_1(X),\lambda_2(X),\ldots,\allowbreak \lambda_n(X))^\top$ and $\delta(X):=(\delta_1(X),\delta_2(X),\ldots,\delta_n(X))^\top$. Specifically for $C \in \mathbb{S}^n_{+}$ (which plays a special role for us), we define $\lambda_{\max} := \lambda_{1}(C)$,  $\lambda_{\min} := \lambda_{n}(C)$, and  $\mu_{\max}$ ($\mu_{\min}$) as the multiplicity of $\lambda_{\max}$ ($\lambda_{\min}$) as an eigenvalue of $C$. Specifically for a diagonal matrix $D \in \mathbb{S}^n$ (which plays a special role for us), we let $d_{\max}:=\delta_1(D)$  and $d_{\min}:=\delta_n(D)$. For a matrix $H := \left(\begin{smallmatrix}
    U & W\\
    Y & Z
\end{smallmatrix}\right)$, we denote the Schur complement of the block $U$ in  $H$ by $H/U$.

\section{On the equivalence of MESP and 0/1 D-Opt}\label{sec:equivalence}
We begin by formally defining $\MESP$, 0/1 $\DOPT$, and some key special cases.
Then we introduce four mappings between instances of these problems, and we work out their important properties.
\medskip

\noindent{\bf MESP.} 
Let $C$ be a symmetric positive-semidefinite matrix with rows/columns
indexed from $N$.
For $0< s < n$,
we define the \emph{maximum-entropy sampling problem}
\begin{align*}
\tag{MESP$(C,s)$}\label{MESP}
\begin{array}{rl}
&\zmespthing(C,s):= \zmespthing(\MESP(C,s))\\
&\quad :=
\max \left\{\ldet \left(C_{S,S} \right)~:~ S\subset N, |S| = s\right\},\\
&\quad =\max \left\{ \ldet\left( C_{S(x),S(x)}\right)~:~ \mathbf{e}^\top x =s,~ x\in\{0,1\}^n\right\},
\end{array}
\end{align*}
where 
 $C_{S,S}$ denotes the principal submatrix indexed by $S$, and $S(x)$ denotes the support of $x$. 
 For feasibility, we assume that $r:=\rank(C)\geq s$. 
 \ref{MESP} was introduced by \cite{SW}; also see \cite{FLbook} and the many references therein. Briefly, in the Gaussian case,  $\ldet (C_{S,S})$ is
 proportional to the ``differential entropy'' (see \cite{Shannon}) of a vector of random variables 
 having covariance matrix $C_{S,S}$\,.  So \ref{MESP} seeks to
 find the ``most informative'' $s$-subvector from an $n$-vector following a joint 
 Gaussian distribution. \ref{MESP} finds application in many areas, e.g. environmental monitoring (see \cite[Chap. 4]{FLbook}). 
In the context of solving \ref{MESP} to optimality by B\&B, 
we encounter subproblem instances of the form 
\begin{equation*}
\ldet (C_{F_1,F_1} ) + \max \left\{\ldet \left(\bar C_{S,S} \right)~:~  S\subset N\setminus F_0 \setminus F_1,~ |S| = s-|F_1|\right\},
\end{equation*}
where $\bar C$ is the Schur complement of $C_{F_1,F_1}$ 
in $C_{N\setminus F_0,N\setminus F_0}$\,. 
We can immediately see that such a \ref{MESP} subproblem instance is 
just another \ref{MESP} problem instance, with optimal objective value
shifted by $\ldet ( C_{F_1,F_1} )$. 

 \medskip
\noindent{\bf 0/1 D-Opt.}  
 Consider the 0/1 $\DOPT$ problem formulated as 
\begin{align*}\label{DOPT}\tag{D-Opt$(A,B,s)$}
\begin{array}{rl}
\textstyle
&\zdoptthing(A,B,s):=
\zdoptthing(\DOPT(A,B,s))\\
&~ :=
\max \left\{ \ldet \left( A_{S\cdot}^\top A_{S\cdot} + B^\top B\right) ~:~ S\subset N, |S| = s
\right\},\\
&~ =\max\!\left\{
\ldet\! \left( A^\top \Diag(x) A \!+\! B^\top B\right)
~\!\!:\!\!~ \mathbf{e}^\top x = s,\, x\!\in\!\{0,1\}^n
\right\},
\end{array}
\end{align*}
where $A \in \mathbb{R}^{n\times m}$  (with no all-zero rows), $B \in \mathbb{R}^{q\times m}$, $q\geq 1$  (with all-zero rows allowed), 
{\scriptsize$\begin{pmatrix}A\\ B \end{pmatrix}$} has full column rank $m$, and $n>s\geq m-\rank(B)$. In the world of statistics, rows of {\scriptsize$\begin{pmatrix}A\\ B \end{pmatrix}$}
are considered to be ``design points'' encoding levels 
of $m$ factors. We are required to choose all rows of $B$, and we wish to choose, additionally, $s$ distinct rows from $A$. 
Here, $B^\top B$ is known as the ``existing Fisher-Information matrix''.
The goal is to minimize the generalized variance for parameter estimates in a least-squares model, based on the chosen design points. 
The assumption that {\scriptsize$\begin{pmatrix}A\\ B \end{pmatrix}$} has full column rank
implies that every feasible solution of \ref{DOPT} yields an ``identifiable'' regression model
(i.e., the associated least-squares problem has a unique solution). 
The 0/1 $\DOPT$ problem and its variants are truly fundamental in the design of experiments; see \cite{puk}. Also, see \cite{PonteFampaLeeMPB} and the many references therein.

An important point on notation: We use $\zmespthing$ for the optimal objective value of a $\MESP$ instance and $\zdoptthing$ for the optimal objective value of a $\DOPT$ instance. In \S\ref{sec:compare-bounds}, we work with many upper bounds, some for $\zmespthing$ and some for $\zdoptthing$\,. We  use $z$ with a subscript that indicates the bound name for upper bounds on $\zmespthing$ and
$\mathfrak{z}$ with a subscript that indicates the bound name for upper bounds on $\zdoptthing$\,. 

Some terminology: When $C$ is positive definite, we say that \ref{MESP} is
\emph{positive definite}.
When $B$ has only zero rows, we say that
\ref{DOPT} is a \emph{pure 0/1 $\DOPT$ problem}, and this is the version often presented
(see \cite{PonteFampaLeeMPB,hugedoptmohit,mohit2020approximation}, e.g.). The distinction between what we refer to as
\ref{DOPT} and the particular case that is the pure 0/1 $\DOPT$ problem
is often glossed over, but in what follows, we need to be precise about this.
Importantly, when carrying out B\&B on a \ref{DOPT} instance
(pure or not), B\&B subproblems are general instances of \ref{DOPT} (not pure).
When $B$ has full column rank $m$ (equivalently, $B^\top B\in\mathbb{S}^m_{++}$),  \ref{DOPT} is known as a \emph{0/1 D-optimal data-fusion problem}; see \cite{li2022d}.  When attacking any \ref{DOPT} with B\&B, subproblems 
can become 0/1 D-optimal data-fusion problems, and this has been exploited (see \cite[the first paragraph of Sec.~7]{PonteFampaLeeMPB}) 
because 0/1 D-optimal data-fusion specific bounds are available; see \cite{PonteFampaLeeMPB}, and the references therein.

We will demonstrate a
very strong relationship between $\MESP$ and 0/1 $\DOPT$.
For models \ref{MESP} and \ref{DOPT},
and a constant $K\in \mathbb{R}$,
we write \ref{MESP}$+K$ and \ref{DOPT}$+K$
to denote precisely the same models but 
with the constant $K$ added to the objective functions.
In this case, we could say that the resulting
models are essentially identical to the given models.
More abstractly and more generally, we are considering 
combinatorial-optimization problems
for which the feasible sets
are points in $\{0,1\}^n$.
If we have combinatorial optimization 
problems $P$ and $Q$ (which could have different models but) with identical 
feasible sets on $\{0,1\}^n$,
and a constant $K\in\mathbb{R}$,
such that every feasible solution has 
objective value in $Q$ exactly $K$ more than in $P$, we write $Q\equiv P+ K$,
and we say that $P$ and $Q$ are 
\emph{completely equivalent}.

We also consider equivalence
under the affine involution
that takes $x$ to  $\mathbf{e}-x$, which we refer to as the \emph{complementation involution.}
We consider  
a combinatorial-optimization problem $Q$
for which the set of feasible 
solutions is the image of the set of feasible solutions of $P$ under the complementation involution.
If we have a constant $K\in\mathbb{R}$,
such that for every feasible solution $\hat{x}$ of $P$, the objective value of $Q$ at $\mathbf{e}-\hat{x}$ is exactly $K$ more than the objective value of $P$ at $\hat{x}$, we write $Q \circeq P + K$, and we say that $P$ and $Q$ are  \emph{complementation equivalent}.

\cite{AFLW_Using} first pointed out (also see \cite[Sec. 1.6]{FLbook}) a ``scaling'' principle for $\MESP$:
that is, for $\gamma>0$, we have  $\MESP(C,s)\equiv\MESP(\gamma C,s) -s \log \gamma$,
and so $\zmespthing(C,s)= \zmespthing(\gamma C,s) -s \log \gamma$. We note that any upper bound for $\MESP(\allowbreak \gamma C,s)$ minus $s \log \gamma$ is also an upper bound for $\MESP(C,s)$, which we call a \emph{scaled upper bound}.
Also, \cite{AFLW_Using} pointed out (also see \cite[Sec. 1.6]{FLbook}) a ``complementing'' principle for $\MESP$:
If $C$ is positive definite, 
then $\MESP(C,s)\circeq\MESP(C^{-1},n-s) + \ldet (C)$, here employing the complementation involution,
and so $\zmespthing(C,s)\allowbreak =\allowbreak \zmespthing(C^{-1},n-s) + \ldet (C)$.
In what follows, we say that 
$\MESP(C,s)$ and $\MESP(C^{-1},n-s)$ are \emph{complementary $\MESP$ instances}.  We note that any upper bound for $\MESP(C^{-1},n-s)$ plus $\ldet (C)$ is also an upper bound for $\MESP(C,s)$, which we call a \emph{complementary upper bound}.
These two principles, scaling and complementing, are referred to later in this section and in \S\ref{sec:compare-bounds}.
 \cite[Thms. 1 and 2]{li2022d}
established that all 0/1 D-optimal data-fusion instances
(i.e., \ref{DOPT} with $B^\top B$ being positive definite) 
are completely equivalent to positive-definite \ref{MESP} instances
(i.e., \ref{MESP} with positive-definite $C$), and conversely.
In what follows, we will establish in full generality that
all \ref{MESP} instances are complementation equivalent to \ref{DOPT} instances, and conversely.

The equivalences that we wish to discuss, are based on four mappings.

\begin{definition}\label{def:mapping}~ 
\begin{itemize}
\item  
We define $\mathcal{M}$, from arbitrary instances of 0/1 $\DOPT$ to instances of $\MESP$, via $\mathcal{M}(\DOPT(A,B,s)):=\MESP(C,n-s)$,
where $\DOPT(A_{n\times m}\,,B_{q\times m}\,,s)$ is an instance of 0/1 
$\DOPT$ and $C:=I_n-A(A^\top A + B^\top B)^{-1}A^\top$.
Note that $\rank(C)$ is 
equal to $n$ minus the multiplicity of 1 as an eigenvalue of 
$A(A^\top A + B^\top B)^{-1}A^\top$; see Proposition \ref{prop:4sure}, below.

\item 
We define $\mathcal{P}$, from arbitrary 0/1 D-optimal data-fusion instances (i.e., instances of 0/1 $\DOPT$ with $B^\top B\in\mathbb{S}^m_{++}$) to positive-definite instances of $\MESP$,
via $\mathcal{P}(\DOPT(A,B,s)):=\MESP(C,s)$,
where $\DOPT(A_{n\times m}\,,B_{q\times m}\,,s)$ is an instance of 
0/1 $\DOPT$, and $C:=I_n+A(B^\top B)^{-1}A^\top$.

\item 
Given an arbitrary $\MESP$ instance $\MESP(C_{n\times n}\,,s)$ 
and a real Schur decomposition $\Phi\Lambda\Phi^\top$ of $C$,
we define $\mathcal{D}(\MESP(C,s);\Phi):=\DOPT(A,B,n-s)$,
where
$A_{n\times n}:=\Phi(I_n-\frac{1}{\lambda_{\max}}\Lambda)^{\scriptscriptstyle 1/2}$, and $B_{n\times n}:=\frac{1}{\sqrt{\lambda_{\max}}}\Lambda^{\scriptscriptstyle 1/2}$. 
Note that $\rank(B)=\rank(C)$, and $\rank(A)=n-\mu_{\max}$\,. Additionally, $A^\top A+B^\top B=I_n$ and $AA^\top=I_n - \frac{1}{\lambda_{\max}}C$, which are facts that we will use in what follows.

\item 
Given a positive-definite $\MESP$ instance $\MESP(C_{n\times n}\,,s)$ and 
$A\in\mathbb{R}^{n\times m}$ such that $AA^\top:=(1/\lambda_{\min})C-I_n$\,, 
we define $\mathcal{F}(\MESP(C,s);A):=\DOPT(\allowbreak A,I_m\,,s)$, which is a 0/1 D-optimal data-fusion instance.  Note that $m\geq \rank(A)=n-\mu_{\min}$\,.
\end{itemize}
\end{definition}

It is clear that the $C$ produced by the map $\mathcal{P}$ is positive definite, and so we have a legitimate $\MESP$ instance from $\mathcal{P}$.
For the $C$ produced by $\mathcal{M}$, we have the following simple result, the first
part of which demonstrates that the $\MESP$ instance produced from $\mathcal{M}$ is legitimate.

\begin{proposition}\label{prop:4sure}
$I_n-A(A^\top A + B^\top B)^{-1}A^\top$ is positive semidefinite and has maximum eigenvalue less than or equal to one.
\end{proposition}

The formula for $C$
in the definition of the mapping $\mathcal{M}$ appeared in \cite[Sec. 3.4.1]{PonteFampaLeeMPB}, but it was mainly discussed there for the case of $B:=0$ (pure 0/1 $\DOPT$)
and also referred to there
(and in \cite{li2022d}) for the case of $B^\top B\in\mathbb{S}^m_{++}$ (i.e., 0/1 D-optimal data-fusion).  
The mapping $\mathcal{P}$ was used in 
\cite[Thm. 2]{li2022d}.
The mapping $\mathcal{D}$
appears to be completely new.
The mapping $\mathcal{F}$ was used in \cite[Thm. 1]{li2022d}.

Notice that both $\mathcal{M}$ and $\mathcal{P}$ map instances of 
0/1 $\DOPT$ to instances of $\MESP$. But  $\mathcal{P}$ operates 
only on 0/1 D-optimal data-fusion instances of 0/1 $\DOPT$ and always produces positive-definite $\MESP$ instances, while 
the mapping $\mathcal{M}$ works on completely general instances of
0/1 $\DOPT$ and can produce singular $\MESP$ instances. Similarly,
$\mathcal{D}$ and $\mathcal{F}$ map instances of 
$\MESP$ to instances of 0/1 $\DOPT$. But  $\mathcal{F}$ operates 
only on positive-definite instances of  $\MESP$ and always produces 0/1 D-optimal data-fusion instances, while 
our new mapping $\mathcal{D}$ works on completely general instances of
$\MESP$ and can produce instances of 0/1 $\DOPT$ that are not 0/1 D-optimal data-fusion instances.

$\mathcal{M}$ and $\mathcal{D}$ map instances that seek to
select $s$ elements from $N$ to instances that select $n-s$ elements
from $N$. This is in contrast to $\mathcal{P}$ and 
$\mathcal{F}$ which map instances that seek to
select $s$ elements from $N$ to instances that select $s$ elements
from $N$. Nonetheless, we will see that $\mathcal{M}$ generalizes $\mathcal{P}$,
and $\mathcal{D}$ generalizes $\mathcal{F}$. 

$\mathcal{D}$ and $\mathcal{F}$ are both defined with respect to 
a factorization. $\mathcal{D}$  uses a  real Schur decomposition $\Phi\Lambda\Phi^\top$ of $C$ (e.g., not unique when 
eigenvalues of $C$ have multiplicity greater than 1),
while $\mathcal{F}$ uses any factorization $AA^\top$ of $(1/\lambda_{\min})C-I_n$\,.
Next, we will demonstrate that for both $\mathcal{D}$ and $\mathcal{F}$, the optimal values of the resulting 0/1 $\DOPT$ instances
do not depend on the chosen factorizations. 

\begin{lemma}\label{lem:IndependentPhi}
     Consider an arbitrary $\MESP$ instance $\MESP(C,s)$
     and two real Schur decompositions $\Phi_i \Lambda \Phi_i^\top$ of $C$, $i=1,2$.
     Let $\DOPT(A_i\,,B,n-s):=\mathcal{D}(\MESP(C,s);\Phi_i)$,
where
$A_i:=\Phi_i(I_n-\frac{1}{\lambda_{\max}}\Lambda)^{\scriptscriptstyle 1/2}$, and $B:=\frac{1}{\sqrt{\lambda_{\max}}}\Lambda^{\scriptscriptstyle 1/2}$. Then: 
\begin{enumerate}
\item[(i)] every $\hat x\in[0,1]^n$  has the same objective value in
$\DOPT(A_1\,,B,n-s)$ and $\DOPT(A_2\,,B,n-s)$.
\item[(ii)]$\zdoptthing(A_1,B,n-s)=\zdoptthing(A_2,B,n-s)$.
\end{enumerate}
\end{lemma}

\begin{proof}
The objective value of $\hat x$ in $\DOPT(A_i\,,B,n-s)$ is
    \begin{align*}
    &\ldet\left( {A}_i^\top \Diag(\hat x) {A}_i + {B}^\top {B}\right)\\
        &\quad = \ldet\left( {A}_i^\top {A}_i +B^\top B   - {A}_i^\top \Diag(\mathbf{e}-\hat x) {A}_i\right) \\ 
                &\quad =\ldet\left( I_n - {A}_i^\top\Diag(\mathbf{e}-\hat x){A}_i\right) \mbox{\qquad (using ${A}_i^\top {A}_i +B^\top B=I_n$)}\\
                    &\quad = \ldet\left( I_n - \Diag(\mathbf{e}-\hat x)^{\scriptscriptstyle 1/2}{A}_i{A}_i^{\top}\Diag(\mathbf{e}-\hat x)^{\scriptscriptstyle 1/2}\right) \\ 
        &\quad = \ldet\left( I_n + \Diag(\mathbf{e}-\hat x)^{\scriptscriptstyle 1/2}((1/\lambda_{\max}) C-I_n)\Diag(\mathbf{e}-\hat x)^{\scriptscriptstyle 1/2} \right).
    \end{align*}
$(i)$ now follows, 
and we can easily see how $(ii)$ follows from $(i)$. 
\qed \end{proof}

\begin{lemma}\label{independentA}
 Consider a positive-definite $\MESP$ instance $\MESP(C,s)$
     and two  factorizations
     $A_i A_i^\top:=(1/\lambda_{\min})C-I_n$\,, with $A_i\in\mathbb{R}^{n\times m_i}$, for $i=1,2$.
     Let $\DOPT(A_i\,,\allowbreak I_{m_i}\,,s):=\mathcal{F}(\MESP(C,s);  A_i)$.
 Then: 
\begin{enumerate}
\item[(i)] every $\hat x\in[0,1]^n$  has the same objective value in
$\DOPT(A_1\,,I_{m_1}\,,s)$ and $\DOPT(\allowbreak A_2\,,I_{m_2}\,,s)$.
\item[(ii)]$\zdoptthing(A_1\,,I_{m_1}\,,s)=\zdoptthing(A_2\,,I_{m_2}\,,s)$.
\end{enumerate}
\end{lemma}

\begin{proof}
The objective value of $\hat x$ in $\DOPT(A_i\,,I_{m_i}\,,s)$ is
    \begin{align*}
    &\ldet\left( {A}_i^\top \Diag(\hat x) {A}_i + I_n\right)\\
                    &\quad = \ldet\left(\Diag(\hat x)^{\scriptscriptstyle 1/2}{A}_i{A}_i^{\top}\Diag(\hat x)^{\scriptscriptstyle 1/2} +I_n\right) \\ 
        &\quad = \ldet\left(\Diag(\hat x)^{\scriptscriptstyle 1/2}((1/\lambda_{\min}) C-I_n)\Diag(\hat x)^{\scriptscriptstyle 1/2} +I_n\right).
    \end{align*}
$(i)$ now follows, 
and we can easily see how $(ii)$ follows from $(i)$. 
\qed \end{proof}

\noindent Next, we mention two results justifying the definitions of $\mathcal{F}$ and 
$\mathcal{P}$. 

\begin{theorem}[\cite{li2022d}]\label{thm:mesp-fu-equiv}
Consider a positive-definite $\MESP$ instance
$\MESP(C,s)$ and any $A$ such that $AA^\top = (1/\lambda_{\min}) C-I_n$\,. We have $\MESP(C,s)$ $\equiv$ \allowbreak $\mathcal{F}(\MESP(C,s);A)$ $+$ $s\log \lambda_{\min}$\,. 
\end{theorem}

\begin{theorem}[\cite{li2022d}]\label{thm:fu-mesp-equiv}
    Consider the 0/1 D-optimal data-fusion instance $\DOPT(\allowbreak A,B,s)$ (i.e., $B^\top B\in\mathbb{S}^m_{++}$).  Then, we have that $\DOPT(A,B,s) \equiv \mathcal{P}(\DOPT(\allowbreak A,B,s))+\ldet(B^\top B)$. 
\end{theorem}

Next, we will see how, given a 0/1 D-optimal data-fusion instance,
the $\MESP$ instances produced by $\mathcal{P}$ and $\mathcal{M}$ are complementary. 

\begin{theorem}\label{thm:relation_m_p}
    Consider the 0/1 D-optimal data-fusion instance $\DOPT(A,B,s)$ (i.e., $B^\top B\in\mathbb{S}^m_{++}$). Let $\MESP(C,s) := \mathcal{P}(\DOPT(A,B,s))$, and let 
    $\MESP(\tilde C,n-s) := \mathcal{M}(\DOPT(A,B,s))$.
    Then $\tilde C= C^{-1}$, and so
    $\mathcal{P}(\DOPT(A,B,s)) \circeq \mathcal{M}(\DOPT(\allowbreak A,B,s)) + \ldet(C)$.
\end{theorem}

\begin{proof}
Considering the complementing principle, we only need to
see that the inverse of $I_n + A(B^\top B)^{-1}A^\top$ is 
$I_n-A(A^\top A + B^\top B)^{-1}A^\top$, which can be easily verified using the 
Woodbury matrix identity.
\qed \end{proof}

Next, we will see that for a pure 0/1 $\DOPT$ instance
$\DOPT(A,0,s)$,
$\mathcal{M}(\DOPT(\allowbreak A,0,s))$ reduces to a construction supplied by 
\cite[first half of Remark 8]{PonteFampaLeeMPB}.

\begin{theorem}\label{thm:aboutpure}
Consider a pure 0/1 $\DOPT$ instance $\DOPT(A_{n\times m}\,,0,s)$; so $A$ must have full column rank $m$. Let $\MESP(C,n-s):=\mathcal{M}(\DOPT(A,0,s))$,
and let $U\Sigma V^\top$ be the compact singular-value decomposition of $A$
(so $U_{n\times m}$\,, $\Sigma_{m\times m}$\,, and $V_{m\times m}$). 
Then $C=I_n-UU^\top$, which has eigenvalues only 0 and 1. 
\end{theorem}

\begin{proof}
Plug $A=U\Sigma V^\top$ and $B=0$ into $C:=I_n-A(A^\top A + B^\top B)^{-1}A^\top$, and observe that $UU^\top$ and $U^\top U= I_m$ have the same 
nonzero eigenvalues.
\qed \end{proof}

We say that \ref{DOPT} is \emph{linearly related} to 
$\DOPT(AM^{-1},BM^{-1},s)$, for every invertible
 $M\in\mathbb{R}^{m\times m}$. It is 
easy to see that linearly related 0/1 $\DOPT$ instances are completely equivalent. 
Next, we will see how
for any 0/1 $\DOPT$ instance produced by $\mathcal{D}$
applied to a positive-definite $\MESP$ instance,
we can choose a 0/1 $\DOPT$ produced by $\mathcal{F}$
applied to the \emph{complementary} $\MESP$ instance,
so that the two 0/1 $\DOPT$ instances are linearly related.  

\begin{theorem}\label{thm:linear}
For a positive-definite $\MESP$ instance $\MESP(C,s)$,
    let $\Phi\Lambda\Phi^\top$ be any
 real Schur decomposition of $C$. 
\vspace{-5pt}
 \begin{enumerate}
     \item \label{DF} Let $\DOPT(A,B,n-s):=\mathcal{D}(\MESP(C,s);\Phi)$, and let $\DOPT(\tilde{A},I_n\,,n-s):=\mathcal{F}(\MESP(C^{-1},n-s);\tilde{A})$, with $\tilde A:=\Phi(\lambda_{\max}\Lambda^{-1}- I_n)^{\scriptscriptstyle 1/2}$. 
\item \label{FD} Let $\DOPT(\tilde A,I_n\,,s):=\mathcal{F}(\MESP(C,s);\tilde A)$, with $\tilde A:=\Phi(\frac{1}{\lambda_{\min}}\Lambda - I_n)^{\scriptscriptstyle 1/2}$  and $\DOPT(A,B,s):=\mathcal{D}(\MESP(C^{-1},n-s);\Phi)$. 
 \end{enumerate}
\vspace{-5pt}
 In both cases described above, we have that 
$\tilde A= A B^{-1}$, and so 
$\DOPT(A,B,n-s)$ and $\DOPT(\tilde{A},I_n\,,n-s)$ are 
linearly related. 
\end{theorem}

\begin{proof}
To verify item \ref{DF}, we note that   $\tilde A\tilde{A}^\top =\Phi(\lambda_{\max}\Lambda^{-1}- I_n)\Phi^\top = \lambda_{\max} C^{-1}-I_n$\,, so $\mathcal{F}$ is correctly applied, and also
$
\textstyle \tilde{A}   
   =\Phi(\lambda_{\max}\Lambda^{-1}-I)^{\scriptscriptstyle 1/2} 
   = \Phi\left(I_n-\frac{1}{\lambda_{\max}}\Lambda\right)^{\scriptscriptstyle 1/2}\allowbreak {\sqrt{\lambda_{\max}}}\Lambda^{\scriptscriptstyle -1/2} 
   =\textstyle A B^{-1}.
$
Item \ref{FD} can be similarly verified.
\qed
\end{proof}

Next, we provide a pair of results 
that justify the definitions of $\mathcal{D}$ and 
$\mathcal{M}$. 

\begin{theorem}\label{thm:equivD}
For an arbitrary $\MESP$ instance $\MESP(C,s)$ 
 and 
any real Schur decomposition $\Phi\Lambda\Phi^\top$ of $C$, we have 
$\MESP(C,s)\circeq\mathcal{D}(\MESP(C,s);\Phi)+s\log \lambda_{\max}$\,.
\end{theorem}

\begin{proof}
    Let  $\DOPT(A,B,n-s):=\mathcal{D}(\MESP(C,s);\Phi)$. 
    Let $\hat S$ be a feasible solution of \ref{MESP}, and let $\hat T:=N\setminus \hat S$, then we have 
    \begin{align*}
       & \ldet(C_{\hat S,\hat S})  = \ldet\left(((1/{\lambda_{\max}}) C)_{\hat S,\hat S}\right) +s\log\lambda_{\max} \\
       &\quad = \ldet\left( (I_n - AA^\top)_{\hat S,\hat S}\right) +s\log\lambda_{\max}\\ 
        &\quad = \ldet\left( I_s - A_{\hat S\cdot}  A_{\hat S\cdot}^\top\right) +s\log\lambda_{\max}
        = \ldet\left( I_n - A_{\hat S\cdot}^\top  A_{\hat S\cdot}\right) +s\log\lambda_{\max}\\ 
        &\quad = \ldet\left(  I_n - (1/\lambda_{\max})\Lambda  - A_{\hat S\cdot}^\top  A_{\hat S\cdot} + (1/\lambda_{\max})\Lambda\right) +s\log\lambda_{\max}\\ 
        &\quad = \ldet\left( A^\top A  - A_{\hat S\cdot}^\top  A_{\hat S\cdot} +  (1/\lambda_{\max})\Lambda\right) +s\log\lambda_{\max}\\ 
        &\quad = \ldet\left(A_{\hat T\cdot}^\top A_{\hat T\cdot} +  B^\top B\right) +s\log\lambda_{\max}\,.
    \end{align*}
The result follows.
\qed
\end{proof}

\begin{theorem}\label{thm:equivM}
For an arbitrary 0/1 $\DOPT$ instance $\DOPT(A,B,s)$, 
we have 
$\DOPT(\allowbreak A,B,s)\circeq\mathcal{M}(\DOPT(A,B,s))+\ldet(A^\top A + B^\top B)$. 
\end{theorem}

\begin{proof}
Let
 $\MESP(C,n-s):=\mathcal{M}(\DOPT(A,B,s))$.   
    Let $\hat S$ be a feasible solution of \ref{DOPT}. Let $\hat T:=N\setminus \hat S$  and  $t := |\hat T|$, then we have 
    \begin{align*}
       & \ldet(A_{\hat S\cdot}^\top A_{\hat S\cdot} + B^\top B)
        = \ldet( A^\top A +  B^\top B -  A_{\hat T\cdot}^\top A_{\hat T\cdot})\\
        &\quad= \ldet( I_m -  ( A^\top A +  B^\top B)^{\scriptscriptstyle -1/2}A_{\hat T\cdot}^\top A_{\hat T\cdot}( A^\top A +  B^\top B)^{\scriptscriptstyle -1/2})+\ldet(A^\top A + B^\top B)\\
        &\quad = \ldet( I_t -  A_{\hat T\cdot}( A^\top A +  B^\top B)^{-1} A_{\hat T\cdot}^\top)+\ldet(A^\top A + B^\top B)\\
        &\quad = \ldet( (I_n -  A( A^\top A +  B^\top B)^{-1} A^\top)_{\hat T,\hat T})+\ldet(A^\top A + B^\top B)\\
        &\quad = \ldet( {C}_{\hat T,\hat T})+\ldet(A^\top A + B^\top B).
    \end{align*}
The result follows.
\qed \end{proof}

$\MESP$ is
NP-hard when $C$ (or $C^{-1}$) has only entries of 
$0,1,3n$ (see \cite{KLQ}).
$\MESP$ is NP-hard when $C$ (or $C^{-1}$) is
an ``arrowhead matrix'' (see \cite{Ohsaka}). 
As a consequence of Theorem \ref{thm:equivM}, we can identify another class of NP-hard $\MESP$ instances.

\begin{corollary}\label{cor:mespnphard}
     $\MESP$ is NP-hard even when the covariance matrix $C$ is a rank-deficient matrix with all positive eigenvalues equal to one.
\end{corollary}
\begin{proof}
    The class of pure $\DOPT$ instances is NP-hard; see \cite{welch1982algorithmic}. 
    From Theorem \ref{thm:equivM}, we can reduce $\DOPT(A,0,n-s)$ to  $\MESP(C,s)$, where $C := I_n - A(A^\top A)^{-1}A^\top$. It is easy to see that $C$ is a rank-deficient matrix  with all positive eigenvalues equal to one; see Theorem \ref{thm:aboutpure}. If $A$ is rational, then
    the reduction is polynomial time.
\qed \end{proof}

In the next two results, we demonstrate that the maps $\mathcal{D}$ and $\mathcal{M}$ 
are kind of inverses of one another.

\begin{theorem}\label{thm:doubleM}
For an arbitrary $\MESP$ instance $\MESP(C,s)$ 
and 
any real Schur decomposition $\Phi\Lambda\Phi^\top$ of $C$, we have that
$\mathcal{M}(\mathcal{D}(\MESP(C,s);\Phi))$ is precisely  $\MESP(\allowbreak (1/\lambda_{\max})C,s)$.
\end{theorem}

\begin{proof}
Let
$
\mbox{D-Opt}(A, B, n-s):= \mathcal{D}(\mbox{MESP}(C,s);\Phi)$. 
By the definition of $\mathcal{D}$, we 
have   
$A:=\Phi(I-\frac{1}{\lambda_{\max}}\Lambda)^{\scriptscriptstyle 1/2}$, $B:=\frac{1}{\scriptstyle\sqrt{\lambda_{\max}}}\Lambda^{\scriptscriptstyle 1/2}$.
Now, let
$
\mbox{MESP}(\tilde C,s):=\mathcal{M}(\mbox{D-Opt}(A, B, n-s)).
$
By the definition of $\mathcal{M}$, we 
have   
\begin{align*}
\tilde{C} :=& I_n-A( A^\top  A +  B^\top  B)^{-1} A^\top\\
=&\textstyle I_n-\Phi(I_n-\frac{1}{\lambda_{\max}}\Lambda)^{\scriptscriptstyle 1/2}
\left((I-\frac{1}{\lambda_{\max}}\Lambda)^{\scriptscriptstyle 1/2}\Phi^\top \Phi (I_n-\frac{1}{\lambda_{\max}}\Lambda)^{\scriptscriptstyle 1/2}
+ \frac{1}{\lambda_{\max}} \Lambda\right)^{-1}\\
&\qquad\qquad\qquad\qquad\qquad\qquad\times
\textstyle\left(I_n-\frac{1}{\lambda_{\max}}\Lambda\right)^{\scriptscriptstyle 1/2}\Phi^\top\\
=&\textstyle I_n -\Phi(I_n-\frac{1}{\lambda_{\max}}\Lambda)\Phi^\top
= \frac{1}{\lambda_{\max}} C.
\end{align*}
The result follows.
\qed \end{proof}

\begin{theorem}\label{thm:doubleD}
For an arbitrary 0/1 $\DOPT$ instance
$\DOPT(A,B,s)$, let
$C := I_n - A(A^\top A + B^\top B)^{-1}A^\top$, and let $\Phi \Lambda \Phi^\top$ be any real Schur decomposition of $C$. Then
$\mathcal{D}(\mathcal{M}(\DOPT(A,B,s));\Phi) +\ldet(A^\top A + B^\top B)
+s\log\lambda_{\max}\equiv \DOPT(A,B,s)$. 
\end{theorem}

\begin{proof}
We can see that $\MESP(C,n-s) := \mathcal{M}(\DOPT(A,B,s))$.
Now, let $\DOPT(\allowbreak \tilde{A},\tilde{B},s) := \mathcal{D}(\MESP(C,n-s);\Phi)$.
 By the definition of $\mathcal{D}$, we have $\tilde{A}:=\Phi(I_n- 
 \frac{1}{\lambda_{\max}}\Lambda)^{\scriptscriptstyle 1/2}$, and $\tilde{B}:= 
 \frac{1}{\sqrt{\lambda_{\max}}}\Lambda^{\scriptscriptstyle 1/2}$. 
Let $\hat S$ be a feasible solution of \ref{DOPT},
and let $\hat T:=N\setminus \hat S$. Then, from the proof for 
Theorem \ref{thm:equivM}, 
we see that 
\begin{align}\label{dopt2mesp}
        \ldet\left(A_{\hat S\cdot}^\top A_{\hat S\cdot} +  B^\top B\right)=\ldet(C_{\hat T,\hat T}) +\ldet(A^\top A + B^\top B).
\end{align}
Also, noting that $\hat T$ is a feasible solution of $\MESP(C,n-s)$, we can see  from the proof for 
Theorem \ref{thm:equivD} 
that 
\begin{align}\label{mesp2dopt}
        \ldet(C_{\hat T,\hat T})=\ldet\left(\tilde{A}_{\hat S\cdot}^\top \tilde{A}_{\hat S\cdot} +  \tilde{B}^\top \tilde{B}\right)  
        + s\log \lambda_{\max}\,.
\end{align}

Finally, the result follows from \eqref{dopt2mesp} and \eqref{mesp2dopt}.
\qed \end{proof}

\begin{remark}
In Theorem \ref{thm:doubleM} we observed that the $\MESP$ instance  $\MESP(\tilde{C},s):=\mathcal{M}(\mathcal{D}(\MESP(C,s);\Phi))$, obtained by applying $\mathcal{D}$ and $\mathcal{M}$ in sequence, is linearly related to $\MESP(C,s)$ (where we consider the real Schur decomposition $\Phi\Lambda\Phi^\top$ of $C$). From Theorem \ref{thm:doubleD}, we note that starting from an 0/1 $\DOPT$ instance and applying $\mathcal{M}$ and $\mathcal{D}$ in sequence, we do not have a similar result. The instances    $\DOPT(\tilde{A},\tilde{B},s)=\mathcal{D}(\mathcal{M}(\DOPT(A,B,s));\Phi)$ and  $\DOPT(A,B,s)$, are \emph{not} linearly related, as  $\tilde{A}\tilde{A}^\top = I_n - \frac{1}{\lambda_{\max}}(I_n - A(A^\top A+ B^\top B)^{-1}A^\top)$ and $\tilde{B}\tilde{B}^\top  = \frac{1}{\lambda_{\max}}(I_n - \Phi^\top A(A^\top A+ B^\top B)^{-1}A^\top\Phi)$ (considering, in this case, the real Schur decomposition $\Phi\Lambda\Phi^\top$ of $C:=I_n-A(A^\top A+B^\top B)^{-1}A^\top$).
\end{remark}

In Figure 
\ref{fig:picture2}, we have depicted some key points 
established in this section. The subfigure in the top
indicates, starting from
a 0/1 $\DOPT$ instance, what $\mathcal{M}$ leads to, and in the particular case of 
a 0/1 D-optimal data-fusion instance, what $\mathcal{P}$ leads to, and the relationship between these $\MESP$ instances created, when we start with a 0/1 D-optimal data-fusion instance.  The subfigure in the bottom
indicates, starting from
a $\MESP$ instance, what $\mathcal{D}$ leads to, and in the particular case of 
a positive-definite $\MESP$ instance, what $\mathcal{F}$ leads to, and the relationship between these 0/1 $\DOPT$ instances created, when we start with a positive-definite $\MESP$ instance, and complement the instance before applying $\mathcal{F}$.  


\begin{figure}[ht!]
\scalebox{0.75}{
\begin{tikzpicture}
\bf\color{black}\small
\centering

 \filldraw[fill=gray!20!white, draw=black]  (0,1) ellipse (3.3cm and 3.0cm);

\filldraw[fill=gray!20!white, draw=black] (10,1) ellipse (3.3cm and 3.0cm);

\filldraw[fill=gray!35!white, draw=black](0,-0.2) ellipse (2.2cm and 1.7cm);

\node[circle,  text=black, minimum size=0.1cm]    (top1)  at ( 0, 3.6 ) {$\underline{0/1\; \DOPT}$};

\node[circle, text=black, minimum size=0.1cm]    (top2)  at ( 10, 3.6) {$\underline{\MESP}$};

\node[circle, text=black, minimum size=0.1cm]    (l3)  at ( 0, 2 ) {$\DOPT\left(A,B,s\right)$};

\node[circle, text=black, minimum size=0.1cm]    (l2)  at ( 0, -0.5 ) {$\DOPT(A,B,s)$};

\node[circle, text=black, minimum size=0.1cm]    (l2a)  at ( 0, -0.9 ) {$B^\top B\in \mathbb{S}^m_{++}$};

\node[circle, text=black, minimum size=0.1cm]    (r3)  at ( 0, 1 ) {$\underline{\rm{Data\; Fusion}}$};

\node[circle, text=black, minimum size=0.1cm]    (r3)  at ( 10, 2 ) {$\MESP(I_n-A(A^\top A+B^\top B)^{-1}A^\top, n-s)$};

\node[circle, text=black, minimum size=0.1cm]    (r2)  at ( 10, -0.5 ) {$\MESP(I_n+A(B^\top B)^{-1}A^\top,s)$};

\node[circle, text=black, minimum size=0.1cm]    (top1)  at ( 5, 2.3 ) {$\mathcal{M}$};

\draw [-{Stealth[length=5mm]}] (l3) -- (r3);

\node[circle, text=black, minimum size=0.1cm]    (top1)  at ( 5, -0.2) {$\mathcal{P}$};

\draw [-{Stealth[length=5mm]}] (l2) -- (r2);

\node[circle, text=black, minimum size=0.1cm]    (top1)  at ( 9, 0.8) {$\rm{complement}$};

\node[circle, text=black, minimum size=0.1cm]    (top1)  at ( 10.9, 0.8) {$\rm{(Thm. ~\ref{thm:relation_m_p})}$};

\draw [ double=gray!20!white, double distance=2pt, -{Stealth[length=3mm]}] (10,0) -- (10,1.7);

\end{tikzpicture}

}
\end{figure}

\vspace{-0.8in}

\begin{figure}[ht!]
\scalebox{0.75}{
\begin{tikzpicture}
\bf\color{black}\small
\centering

\filldraw[fill=gray!20!white, draw=black] (0,1) ellipse (3.3cm and 3.0cm);

\filldraw[fill=gray!20!white, draw=black] (10,1) ellipse (3.3cm and 3.0cm);

\filldraw[fill=gray!35!white, draw=black] (10,-0.2) ellipse (2.2cm and 1.7cm);

\node[circle, text=black, minimum size=0.1cm]    (top1)  at ( 0, 3.6) {$\underline{0/1\; \DOPT}$};
\node[circle, text=black, minimum size=0.1cm]    (top2)  at ( 10, 3.6) {$\underline{\MESP}$};

\node[circle, text=black, minimum size=0.1cm]    (l3)  at ( 0, 2 ) {$\DOPT\left(\Phi(I_n-\frac{1}{\lambda_{\max}}\Lambda)^\frac{1}{2},\frac{1}{\sqrt{\lambda_{\max}}}\Lambda^{\frac{1}{2}},n-s\right)$};

\node[circle, text=black, minimum size=0.1cm]    (l2)  at ( 0, 0.2 ) {$\DOPT(A,I_n,n-s)$};

\node[circle, text=black, minimum size=0.1cm]    (l1)  at ( 0, -1.5 ) {$\DOPT(A,I_n,s)$};

\node[circle, text=black, minimum size=0.1cm]    (r3)  at ( 10, 1 ) {$\underline{\rm{pos\mbox{-}def }\;\MESP}$};

\node[circle, text=black, minimum size=0.1cm]    (r3)  at ( 10, 2 ) {$\MESP(C,s)$};

\node[circle, text=black, minimum size=0.1cm]    (r2)  at ( 10, 0.2 ) {$\MESP(C^{-1},n-s)$};

\node[circle, text=black, minimum size=0.1cm]    (r1)  at ( 10, -1.5 ) {$\MESP(C,s)$};

\node[circle, text=black, minimum size=0.1cm]    (top1)  at ( 5, 2.3 ) {$\mathcal{D}$};

\node[circle, text=black, minimum size=0.1cm]    (top1)  at ( 5, 1.7 ) {$\Phi\Lambda\Phi^\top:=C$};

\draw [-{Stealth[length=5mm]}] (r3) -- (l3);

\node[circle, text=black, minimum size=0.1cm]    (top1)  at ( 5, 0.5 ) {$\mathcal{F}$};

\node[circle, text=black, minimum size=0.1cm]    (top1)  at ( 5, -0.1 ) {$A:=\Phi(\lambda_{\max}\Lambda^{-1}- I_n)^{\scriptscriptstyle 1/2}$}; 

\node[circle, text=black, minimum size=0.1cm]    (top1)  at ( 5, -0.5 ) {$\Phi\Lambda\Phi^\top:=C$};

\draw [-{Stealth[length=5mm]}] (r2) -- (l2);

\node[circle, text=black, minimum size=0.1cm]    (top1)  at ( 5, -1.2 ) {$\mathcal{F}$};

\node[circle, text=black, minimum size=0.1cm]    (top1)  at ( 5, -1.8 ) {$AA^\top:=\frac{1}{\lambda_{\min}}C-I_n$};

\draw [-{Stealth[length=5mm]}] (r1) -- (l1);

\draw [double=gray!35!white, double distance=2pt, -{Stealth[length=3mm]}] (10,-1.2) -- (10,0);

\node[circle, text=black, minimum size=0.1cm]    (top1)  at ( 9, -0.7 ) {$\rm{complement}$};

\draw [double=gray!20!white, double distance=2pt, -{Stealth[length=3mm]}] (0,0.5) -- (0,1.5);

\node[circle, text=black, minimum size=0.1cm]    (top1)  at ( 0.95, 0.9 ) {$\rm{(Thm. ~\ref{thm:linear})}$};

\node[circle, text=black, minimum size=0.1cm]    (top1)  at ( -0.7, 1.1 ) {$\rm{linearly}$};

\node[circle, text=black, minimum size=0.1cm]    (top1)  at ( -0.7, 0.8 ) {$\rm{related}$};

\end{tikzpicture}
}
\caption{From 0/1 D-Opt and from MESP}\label{fig:picture2}
\end{figure}


\section{Transferring upper bounds for MESP and 0/1 D-Opt}\label{sec:compare-bounds}

By applying any of the many basic upper bounds on the optimal value of a $\MESP$ instance to $\mathcal{M}(\DOPT(A,B,s))$, we can obtain many upper bounds for the optimal value of 
$\DOPT(A,B,s)$,  via Theorem \ref{thm:equivM}. On the other side, there are  only a few basic upper bounds for the optimal value of a $\DOPT$ instance, so by applying them to $\mathcal{D}(\MESP(C,s))$,
we obtain only a few upper bounds for the optimal value of $\MESP(C,s)$, via Theorem \ref{thm:equivD}. In what follows, we describe the
basic upper bounds, and then we develop their properties
in the context of  $\mathcal{M}$, $\mathcal{P}$, $\mathcal{D}$ and $\mathcal{F}$.


\subsection{Basic MESP upper bounds}\label{sec:MESPbounds}

We consider several basic upper bounds from the literature for \ref{MESP}:
\begin{itemize}
    \item \textbf{NLP bound.} Consider given parameters $d,p \in \mathbb{R}^n_{++}$ and $\gamma > 0$. Let $D := \Diag(d)$. Then, the (scaled) \emph{NLP bound}  proposed in \cite{AFLW_IPCO,AFLW_Using}  is given by 
\begin{align}\label{nlp_bound}\tag{NLP}
&\hypertarget{znlptargetmesp}{\znlpthingmesp}(C,s) := \max \left\{\vphantom{\Big)}\ldet\big(\Diag((\gamma d)^x) + \gamma \Diag(x^{p/2})(C-D)\Diag(x^{p/2}) \big)\right. \\[-13pt] 
&\quad \left. - s\log(\gamma) ~~
  \mathbf{e}^\top x\!=\!s,~ 
x\in[0,1]^n\vphantom{\Big)}\right\},\nonumber 
\end{align}

\begin{theorem}[\protect{\cite[Cor. 3.5.9]{FLbook}}]\label{cor:best-p-nlp}
    Let $D\succeq C$ and $\gamma>0$. The minimum upper bound for \ref{MESP}, given by $\znlpmesp$\,,
     is obtained by 
    setting, for $i \in N$:
    \vspace{-5pt}
    \begin{align*}
        p_i:=\begin{cases}
        1, &    \text{for }\gamma d_i\leq 1;
        \\
        \left(1+\sqrt{1+4\log(\gamma d_i)}\right)^2/4,\quad &\text{for }\gamma d_i> 1.
    \end{cases}
    \end{align*}
\end{theorem}

\cite{AFLW_Using} proposed the following three strategies for selecting the parameters in \ref{nlp_bound}, aiming at obtaining the best possible bound for \ref{MESP} while assuring convexity of \ref{nlp_bound}. 
Then we determine $p$ using Theorem \ref{cor:best-p-nlp}.

\noindent --    \hypertarget{Identitytarget}{\Identity}:
 $  
        D:= \rho I,~ \rho:=\lambda_{\max}\,,~ \gamma:= 1/\rho.
$   
By Theorem \ref{cor:best-p-nlp}, we have $p=
\mathbf{e}$.

\noindent --  
    \hypertarget{Diagonaltarget}{\Diagonal}: 
$ 
        D\!:= \!\rho\Diag(C),\rho\!:=\!\lambda_{1}(\Diag(C)^{\scriptscriptstyle -1/2}C\Diag(C)^{\scriptscriptstyle -1/2}),\gamma\in [1/d_{\max},\allowbreak 1/d_{\min}].
$ 

    \par
\noindent --  
\hypertarget{Tracetarget}{\Trace}: 
\vspace{-20.5pt}
\begin{equation}\label{eq:nlp-tr} \quad\qquad\,
        D\!:=\! \argmin\{\Tr(Y):Y-C\succeq 0,\,Y \text{ diagonal}\},\gamma\!\in\! [1/d_{\max},1/d_{\min}].
    \end{equation}

We define as \hypertarget{NLPIdtarget}{\NLPIdsym}, \hypertarget{NLPDitarget}{\NLPDisym} and \hypertarget{NLPTrtarget}{\NLPTrsym}, the special cases of \ref{nlp_bound} where the parameters $d$ and $\gamma$ are defined respectively by the
\Identity, \Diagonal, or \Trace~strategy, and $p$ is defined by Theorem \ref{cor:best-p-nlp}. We denote the optimal objective value of the three problems, respectively by 
$\hypertarget{znlpidtargetmesp}{\znlpidthingmesp}$\,, $\hypertarget{znlpditargetmesp}{\znlpdithingmesp}$\,, $\hypertarget{znlptrtargetmesp}{\znlptrthingmesp}$\,.

\item \textbf{Factorization bound.}
Next, we consider the ``factorization bound''.
It is based on
 the following fundamental lemma. 

\begin{lemma}[\protect{\cite[Lem. 13]{Nikolov}}]\label{Ni13}
 Let $\omega\in\mathbb{R}_+^k$ satisfy $\omega_1\geq \omega_2\geq \cdots\geq \omega_k$\,, define $\omega_0:=+\infty$, and let $s$ be an integer satisfying
 $0<s\leq k$. Then there exists a unique integer $i$, with $0\leq i< s$, such that
$ \omega_{i }>\textstyle\frac{1}{s-i }\textstyle\sum_{\ell=i+1}^k \omega_{\ell}\geq \omega_{i+1}\,.$
\end{lemma}

Suppose that  $\omega\in\mathbb{R}^k_+$ with  
$\omega_1\geq\omega_2\geq\cdots\geq\omega_k$\,.
Let $\hat\imath$ be the unique integer defined by Lemma \ref{Ni13}. We define
\begin{equation*} 
\phi_s(\omega):=\textstyle\sum_{\ell=1}^{\hat\imath} \log \omega_\ell + (s - \hat\imath)\log\left(\frac{1}{s-{\hat\imath}} \sum_{\ell=\hat\imath+1}^{k}
\omega_\ell\right).
\end{equation*}
For $X\in\mathbb{S}_{+}^k$~, we define the \emph{$\Gamma$-function} as
\begin{equation*} 
\Gamma_s(X):= \phi_s(\lambda(X)),
\end{equation*}
and for $X\in\mathbb{S}^k_+$ and $\kappa > 0$, we define the \emph{$\Gamma^+$-function} as
\begin{equation}\label{def:gamma+}
\Gamma_s^+(X):= \phi_s(\lambda(X) + \kappa \,
\mathbb{I}_s).
\end{equation}
where $\mathbb{I}_s\in\mathbb{R}^k$ has the first $s$ elements equal to 1 and the others equal to 0. 
For  $C\in\mathbb{S}^n_{+}$\,, we factorize $C=FF^\top$,
with $F\in \mathbb{R}^{n\times p}$, for some $p$ satisfying $r\le p \le n$. 
The \emph{factorization bound} from 
\cite{Nikolov} (also see \cite{Weijun,FLbook,FactPaper})
is
\begin{align*}\label{ddfact}\tag{\mbox{DDFact}}
&\textstyle \hypertarget{zgammatarget}
\textstyle \hypertarget{zgammatarget}{\zgammathing} (C,s):=\hypertarget{zgammatarget}{\zgammathing} (C,s;F)\\
&\quad :=\max \left\{ \Gamma_s(F^\top \Diag(x)F) \, : \, \mathbf{e}^\top x=s,~
x\in[0,1]^n
\right\}, 
\end{align*}
where we can omit $F$ in notation for $\hypertarget{zgammatarget}{\zgammathing}(C,s;F)$ because the factorization bound does not depend on which factorization of $C$ is used (see \cite{FactPaper}).

For  $C\in\mathbb{S}^n_{++}$\,, we factorize  $C - \lambda_{\min} I_n= GG^\top$,
with  $G \in \mathbb{R}^{n\times q}$, for some $q$ satisfying $\rank(G)\!\leq\! q \! \leq \! n$, where $\rank(G)$ is $n$ minus the multiplicity of the minimum eigenvalue of  $C$. 
The \emph{augmented factorization bound} from \cite{li2025augmented} is
\begin{align*}\label{ddfact_plus}\tag{\mbox{DDFact$_+$}}
&\textstyle \hypertarget{zgammaplustarget}
\textstyle \hypertarget{zgammaplustarget}{\zgammaplusthing}(C,s):=\hypertarget{zgammaplustarget}{\zgammaplusthing}(C,s;G)\\
&\quad :=\max \left\{ \Gamma^+_s(G^\top \Diag(x)G) \, : \, \mathbf{e}^\top x=s,~
x\in[0,1]^n
\right\}, 
\end{align*}
where we set $\kappa:=\lambda_{\min}$ in \eqref{def:gamma+}. We omit the parameter $G$ in the simplified notation for $\hypertarget{zgammaplustarget}{\zgammaplusthing}(C,s;G)$ because, also based on a result of \cite{FactPaper}, the augmented factorization bound does not depend on the chosen factorization of $C-\lambda_{\min}I_n$\,.

\item \textbf{BQP bound.} Given $\gamma>0$. The (scaled) \emph{BQP bound} proposed in  \cite{Anstreicher_BQP_entropy} (also see \cite{FLbook}) is given by  
\begin{align}\label{bqp_original}\tag{BQP}
\textstyle
\hypertarget{zbqptarget}
\textstyle
\hypertarget{zbqptarget}{\zbqpthing}(C,s) :=\max \{\ldet&\left(\gamma C\circ X + \Diag(\mathbf{e}-x) \right) - s\log(\gamma) \, :\\& \, \mathbf{e}^\top x\!=\!s,\,X\mathbf{e}\!=\!sx,\,x\!=\!\diag(X),\,X\!\succeq\!xx^\top\}.\nonumber
\end{align}

\item 
\textbf{linx bound.} Given $\gamma >0$.
The (scaled) \emph{linx bound} proposed in  \cite{Kurt_linx} (also see \cite{FLbook,FactPaper}) is given by  
\begin{align*}\label{linx}\tag{linx}
\textstyle
&
\hypertarget{zlinxtarget}{\zlinx}(C,s) :=\max \left\{\textstyle \frac{1}{2}(\ldet\left(\gamma C\Diag(x) C + \Diag(\mathbf{e}\!-\!x) \right) -s\log(\gamma))\right.\\
&\qquad\qquad\qquad\qquad\qquad  \left.~:~ \mathbf{e}^\top x\!=\!s,~
x\in[0,1]^n
\right\}.
\end{align*}

\item \textbf{Spectral bound.} The \emph{spectral bound} proposed in  \cite{KLQ} is given by 
\begin{align*}\label{spectral-mesp}\tag{\mbox{$\mathcal{S}(\MESP)$}}
\textstyle
&
\textstyle\hypertarget{zspectraltarget}{\zspectral}(C,s) :=\sum_{\ell = 1}^s \log \lambda_{\ell}(C).
\end{align*}

\item \textbf{Diagonal bound.} The \emph{diagonal bound} proposed in  \cite{HLW} is given by 
\begin{align*}\label{hadamard-mesp}\tag{\mbox{diag}}
\textstyle
&
\textstyle\hypertarget{zdiagonaltarget}
{\zdiagonal}(C,s) := \sum_{\ell = 1}^s \log \delta_{\ell}(C).
\end{align*}

\end{itemize}


\subsection{Basic 0/1 D-Opt upper bounds}\label{sec:doptbounds}

We consider the few basic bounds from the literature for  \ref{DOPT}:

\begin{itemize}
\item \textbf{Natural bound.} The most obvious upper bound for \ref{DOPT}
is the \emph{natural bound} from \cite{Welch}; also see \cite{PonteFampaLeeMPB}.
\begin{align*}\label{natural_bound}\tag{\mbox{$\mathcal{N}$}} 
\textstyle \hypertarget{znaturaltarget} 
\textstyle \hypertarget{znaturaltarget}{\znaturalthing}(A,B,s):= 
&\max\! \left\{ \ldet \left(A^\top \Diag(x) A + B^\top B \right) ~\!\!:\!\!~ \mathbf{e}^\top x\!=\!s,~\!\! 
x\in[0,1]^n
\right\}.
\end{align*}

\item \textbf{Spectral bound.} For 0/1 D-optimal data-fusion instances \ref{DOPT} (i.e.,  with $B^\top B\in\mathbb{S}^m_{++}$),
the \emph{spectral bound} from  \cite{KoLeeWayne}, is defined as
\begin{equation}\label{spectral}\tag{\mbox{$\mathcal{S}(\DOPT)$}}
       \textstyle  \hypertarget{zspectraldopttarget}
       \textstyle  \hypertarget{zspectraldopttarget}{\zspectraldoptthing}(A,B,s):= \ldet(B^\top B) +  \sum_{i=1}^{s} \log \lambda_i\left(I_n + A (B^\top B)^{-1}A^\top\right).
    \end{equation}

\item \textbf{Hadamard bound.}   For 0/1 D-optimal data-fusion instances \ref{DOPT} (i.e.,  with $B^\top B\in\mathbb{S}^m_{++}$),
\cite{KoLeeWayne} (also see \cite{Welch})
introduced a \emph{Hadamard bound}, which \cite[Lem. 20]{PonteFampaLeeMPB} reformulate as
\begin{equation}\label{hadamard}\tag{\mbox{$\mathcal{H}$}}
       \textstyle  \hypertarget{zhadamarddopttarget}
       \textstyle  \hypertarget{zhadamarddopttarget}{\zhadamarddoptthing}(A,B,s):=\ldet(B^\top B) +  \sum_{i=1}^{s} \log \delta_i\left(I_n + A (B^\top B)^{-1}A^\top\right).
    \end{equation}

\end{itemize}


\subsection{Transferring upper bounds}\label{subsec:tranf}~

From Theorem \ref{thm:equivM}, we can see that 
given a 0/1 $\DOPT$ instance $\DOPT(A,B,s)$  
with optimal value $\zdoptthing(A,B,s)$, and the $\MESP$ instance $\MESP(C,n-s):= \mathcal{M}(\DOPT(A,B,s))$ with optimal value $\zmespthing(C,n-s)$, then any upper bound for $\zmespthing(C,n-s)$ plus $\ldet(A^\top A+B^\top B)$ is an upper bound for $\zdoptthing(A,B,s)$. 
Similarly, from Theorem \ref{thm:fu-mesp-equiv}, we can see that 
given a 0/1 D-optimal data-fusion  $\DOPT(A,B,s)$ (i.e., $B^\top B\in\mathbb{S}^n_{++}$),  
with optimal value $\zdoptthing(A,B,s)$, and the $\MESP$ instance $\MESP(C,s):= \mathcal{P}(\DOPT(A,B,s))$ with optimal value $\zmespthing(C,s)$, then any upper bound for $\zmespthing(C,s)$ plus $\ldet(B^\top B)$ is an upper bound for $\zdoptthing(A,B,s)$.
But we can see, via Theorem \ref{thm:relation_m_p}, that any bound that can be achieved
by applying a basic MESP upper bound to $\mathcal{P}(\DOPT(A,B,s))$, for a 0/1 D-optimal data-fusion instance, can equally be achieved by applying the same MESP bound to the complement of 
$\mathcal{M}(\DOPT(A,B,s))$. Therefore, in what follows, it suffices to 
study the application of  basic MESP bounds to $\mathcal{M}(\DOPT(A,B,s))$, which
applies to general 0/1 D-Opt, and we do not need to study $\mathcal{P}$ in this context. 

From Theorem \ref{thm:equivD}, we can see that 
given a $\MESP$ instance $\MESP(C,s)$
with optimal value $\zmespthing(C,s)$, 
any real Schur decomposition $\Phi\Lambda\Phi^\top$ of $C$, and the  0/1 $\DOPT$ instance $\DOPT(A,B,n-s):=\mathcal{D}(\MESP(C,s);\Phi)$  
with optimal value $\zdoptthing(A,B,n-s)$, then any upper bound for $\zdoptthing(A,B,n-s)$ plus $s\log \lambda_{\max}$ is an upper bound for $\zmespthing(C,s)$. Similarly, 
from Theorem \ref{thm:mesp-fu-equiv}, we can see that 
given a positive-definite $\MESP$ instance $\MESP(C,s)$
with optimal value $\zmespthing(C,s)$, 
any factorization $AA^\top:=\frac{1}{\lambda_{\min}}C-I_n$\,, and the  0/1 $\DOPT$ instance $\DOPT(A,I_n\,,s):=\mathcal{F}(\MESP(C,s);A)$  
with optimal value $\zdoptthing(A,I_n\,,s)$, then any upper bound for $\zdoptthing(A,I_n\,,s)$ plus $s\log \lambda_{\min}$ is an upper bound for $\zmespthing(C,s)$. 
But we can see, via Theorem \ref{thm:linear} (Item 2),
 that for applying a basic 0/1 D-Opt upper bound to a D-Opt instanced derived from a positive-definite MESP instance, we can confine our attention to $\mathcal{D}$, and we do not need to study $\mathcal{F}$ in this context, as long as the basic 0/1 D-Opt upper bound
shifts by exactly the amount that the optimal value shifts. 
Recall that if $M$ is invertible, then
\ref{DOPT} is \emph{linearly related} to 
$\DOPT(AM^{-1},BM^{-1},s)$. It is easy to check that the
$\DOPT(A,B,s) = \DOPT(AM^{-1},BM^{-1},s) + 2\log(|\det(M)|)$,
and that all three of our basic 0/1 D-Opt upper bounds
likewise shift by exactly $2\log(|\det(M)|)$.

From the discussions above, we highlight two  identities 
 that will be used 
 to show the validity of the bounds that we investigate and compare. 

\begin{align}
&\zdoptthing(A,B,s)= \zmespthing(\mathcal{M}(\DOPT(A,B,s)))+\ldet(A^\top A + B^\top B), \tag{\mbox{$Z.\mathcal{M}$}}\label{ZM}\\
&\zmespthing(C,s)=\zdoptthing(\mathcal{D}(\MESP(C,s);\Phi))+s\log \lambda_{\max}\,.\tag{\mbox{$Z.\mathcal{D}$}}\label{ZD}
\end{align}

As we indicated in \S\ref{sec:MESPbounds}, there are six basic 
upper bounds for $\MESP$, and using the identity \ref{ZM}, we 
induce six $\MESP$-based upper bounds for any 0/1 $\DOPT$ instance \ref{DOPT}. In the following, we refer to an upper bound for $\zdoptthing(A,B,s)$ obtained by an upper bound for  $\zmespthing(\mathcal{M}(\DOPT(A,B,s)))$ plus $\ldet(A^\top A + B^\top B)$ as an \emph{$\mathcal{M}$-induced upper bound}. e.g., the $\mathcal{M}$-induced \ref{nlp_bound} bound for $\DOPT(A,B,s)$ is defined as  $\znlpthingmesp(C,n-s) + \ldet(A^\top A + B^\top B)$, where $\MESP(C,n-s):=\mathcal{M}(\DOPT(A,B,s))$.
We briefly summarize facts about the use of these induced $\MESP$-based upper bounds for 0/1 $\DOPT$, explaining what was
previously known and what we will establish below. 

\begin{itemize}
\item NLP. The $\mathcal{M}$-induced \ref{nlp_bound} bound has never directly been described and experimented
with before. In 
Theorem \ref{thm:natural-nlp-equiv}, we establish that for the
case of $\NLPId$, this upper bound 
strictly
dominates the natural bound for \ref{DOPT} when $A$ has full row rank, and
is equal to the natural bound for \ref{DOPT} otherwise. 
Additionally, in Corollary \ref{cor:jacobi}, we give checkable sufficient conditions under which $\NLPTr$ is equal to $\NLPId$, for the $\MESP$ instances that arise from
\emph{pure} 0/1 $\DOPT$ instances, under $\mathcal{M}$.
\item Factorization. The $\mathcal{M}$-induced factorization bound was heavily worked with in  \cite{PonteFampaLeeMPB} (where it is referred to as the ``$\Gamma$-bound for 0/1 D-optimality'').
\item BQP. The $\mathcal{M}$-induced BQP bound has never been described before.
\item linx. The $\mathcal{M}$-induced linx bound was heavily worked with in \cite{PonteFampaLeeMPB}.
\item Spectral. 
In Theorem \ref{thm:spectral_mesp_eq_dopt},
we establish that
for the case in which \ref{DOPT} is a 0/1 D-optimal data-fusion instance
(i.e., $B^\top B\in\mathbb{S}^n_{++}$), 
the $\mathcal{M}$-induced spectral bound is the same as the spectral bound for 0/1 D-optimal data-fusion. For general instances \ref{DOPT}, this upper bound has not been considered before; see 
Remark \ref{rem:generalspec}.
\item Diagonal. In Theorem \ref{thm:diag_comp_to_had},
we establish that
for the case in which \ref{DOPT} is a 0/1 D-optimal data-fusion instance
(i.e., $B^\top B\in\mathbb{S}^n_{++}$),
the  $\mathcal{M}$-induced complementary diagonal bound
is the same as the Hadamard bound for \ref{DOPT}.
The $\mathcal{M}$-induced diagonal bound, which applies to general instances of \ref{DOPT}, has never been described before.
\end{itemize}

As we indicated in \S\ref{sec:doptbounds}, there are three basic 
upper bounds for $\DOPT$, and using the identity \ref{ZD}, we 
induce three 0/1 $\DOPT$-based upper bounds for any $\MESP$ instance \ref{MESP}. Considering the real Schur decomposition $\Phi\Lambda\Phi^\top$ of $C$, we will refer to an upper bound for $\zmespthing(C,s)$ obtained by an upper bound for  $\zdoptthing(\mathcal{D}(\MESP(C,s);\Phi))$ plus $s\log \lambda_{\max}$ as a \emph{$\mathcal{D}$-induced upper bound}. e.g., the $\mathcal{D}$-induced natural bound for $\MESP(C,s)$ is defined as  $\znaturalthing(A,B,n-s) + s\log \lambda_{\max}$, where $\DOPT(A,B,n-s):=\mathcal{D}(\MESP(C,s);\Phi)$.
Below, we briefly summarize facts about these bounds, explaining what was
known and what we will establish.

\begin{itemize}
\item Natural. In Theorem \ref{thm:nlp-natural-equiv}, we establish that  the $\mathcal{D}$-induced natural bound is equal to the $\NLPId$ bound for \ref{MESP}.
\item Spectral. In Theorem \ref{thm:zsmesp_eq_zdopt_pd}, we
establish that for positive-definite \ref{MESP}, the $\mathcal{D}$-induced spectral bound is equal to the spectral bound for \ref{MESP}. 
\item Hadamard. In Theorem \ref{thm:had}, we establish that 
for the case of positive-definite \ref{MESP},  the $\mathcal{D}$-induced Hadamard bound is equal to the complementary diagonal bound for
 \ref{MESP}.  
\end{itemize}
In the following four subsections, we establish the results mentioned above.


\subsubsection{Spectral bounds.} 
 In this subsection, we 
establish close connections between the spectral bound for $\MESP$ and the
spectral bound for 0/1 $\DOPT$.

\begin{theorem}\label{thm:spectral_mesp_eq_dopt}
    Consider the 0/1 D-optimal data-fusion instance $\DOPT(A,B,s)$ (i.e., $B^\top B\in\mathbb{S}^m_{++}$). Let $\MESP(C,n-s):=\mathcal{M}(\DOPT(A,B,s))$. Then,  the $\mathcal{M}$-induced spectral bound is equal to the spectral bound for $\DOPT(A,B,s)$; that is,
    $\zspectral(C,n-s) + \ldet(A^\top A + B^\top B)=\zspectraldopt(A,B,s)$. 
\end{theorem}

\begin{proof} 
    Note that
\begin{align*}
   &\zspectral(C,n-s) = \textstyle\sum_{i=1}^{n-s} \log \lambda_i\left(C\right) =  \ldet(C) - \sum_{i=n-s+1}^{n} \log \lambda_i\left(C\right) \\
    &\qquad = \textstyle \ldet(C) + \sum_{i=1}^{s} \log \lambda_i\left(C^{-1}\right).
\end{align*}
From the Woodbury identity, we have $C^{-1} = I_n + A (B^\top B)^{-1}A^\top$. Next, note that 
\begin{align}
    &\ldet(A^\top A + B^\top B) + \ldet(C)\label{wood}\\
    &~= \ldet(A^\top A + B^\top B)+ \ldet(I_n - A (A^\top A + B^\top B)^{-1}A^\top)\nonumber\\
    &~= \ldet(A^\top A + B^\top B)+ \ldet(I_n - (A^\top A + B^\top B)^{\scriptscriptstyle -1/2}A^\top A(A^\top A + B^\top B)^{\scriptscriptstyle -1/2})\nonumber\\
    &~=  \ldet(B^\top B).\nonumber
\end{align}
The result follows.
\qed \end{proof}

\begin{remark}\label{rem:generalspec}
   From Theorem \ref{thm:spectral_mesp_eq_dopt}, we note that with the mapping $\mathcal{M}$, we generalize the spectral bound for 0/1 D-optimal data-fusion presented in \cite{KoLeeWayne} to all 0/1 D-Opt instances.
\end{remark}

\begin{theorem}\label{thm:zsmesp_eq_zdopt_pd}
    Consider a positive-definite $\MESP$ instance $\MESP(C,s)$ and 
any real Schur decomposition $\Phi\Lambda\Phi^\top$ of $C$. Let $\DOPT(A,B,n-s):=\mathcal{D}(\MESP(\allowbreak C,s);\Phi)$. Then,  the $\mathcal{D}$-induced spectral bound is equal to the spectral bound for $\MESP(C,s)$; that is,
$\zspectraldopt(A,B,n-s)  + s\log\lambda_{\max}= \zspectral(C,s)$. 
\end{theorem}
\begin{proof}
 By 
 the definition of $\mathcal{D}$, we 
have   
$A:=\Phi(I-\frac{1}{\lambda_{\max}}\Lambda)^{\scriptscriptstyle 1/2}$, $B:=\frac{1}{\scriptstyle\sqrt{\lambda_{\max}}}\Lambda^{\scriptscriptstyle 1/2}$.  Note that $\textstyle{A}  ({B}^\top {B})^{-1} {A}^\top = \lambda_{\max}{C}^{-1} - I_n$ and $\ldet({B}^\top {B}) = \ldet((1/\lambda_{\max}) C)$. We have
\begin{align*}
        &\zspectraldopt(A,B,n-s) + s\log \lambda_{\max} \\
        &\quad = \textstyle \sum_{i=1}^{n-s} \log \lambda_i( \lambda_{\max} C^{-1}) +  \ldet((1/\lambda_{\max})C) + s\log\lambda_{\max}\\
        &\quad=  \textstyle \sum_{i=1}^{n-s} \log \lambda_i(C^{-1}) +  \ldet(C)  = \sum_{i=1}^{s} \log \lambda_i(C) = \zspectral(C,s).
\end{align*}
\qed \end{proof}


\subsubsection{Diagonal and Hadamard bounds.} 
 In this subsection, we 
establish a close connection between the diagonal bound for $\MESP$ and the
Hadamard bound for 0/1 $\DOPT$.

\begin{theorem}\label{thm:diag_comp_to_had}
     Consider the 0/1 D-optimal data-fusion instance $\DOPT(A,\allowbreak B,s)$ (i.e., $B^\top B\in\mathbb{S}^m_{++}$). Let $\MESP(C,n-s):=\mathcal{M}(\DOPT(A,B,s))$. Then, the $\mathcal{M}$-induced  complementary diagonal bound is equal to the Hadamard bound for $\DOPT(\allowbreak A,B,s)$; that is,
     $\zdiagonal(C^{-1},s) + \ldet(C) + \ldet(A^\top A + B^\top B)=\zhadamarddopt(A,B,s)$.
\end{theorem}

\begin{proof}
From the Woodbury identity, we have $C^{-1} = I_n + A (B^\top B)^{-1}A^\top$, and via \eqref{wood},
we get 
$\ldet(C) + \ldet(A^\top A + B^\top B)\!=\! \ldet(B^\top B)$. The result follows. \qed \end{proof}

\begin{theorem}\label{thm:had}
    Consider a positive-definite $\MESP$ instance $\MESP(C,s)$ and 
any real Schur decomposition $\Phi\Lambda\Phi^\top$ of $C$. Let $\DOPT(A,B,n-s)$ $:=$ $\mathcal{D}(\MESP(\allowbreak C,s);\Phi)$. Then,  the $\mathcal{D}$-induced  Hadamard bound is equal to the complementary diagonal bound  for $\MESP(C,s)$; i.e.,
$\zhadamarddopt(A,B,n-s) + s\log\lambda_{\max} = \zdiagonal(C^{-1},n-s)~+~\ldet(C)$. 
\end{theorem}

\begin{proof}
 By the definition of $\mathcal{D}$, we have 
$A:=\Phi(I_n-\frac{1}{\lambda_{\max}}\Lambda)^{\scriptscriptstyle 1/2}$ and $B:=\frac{1}{\sqrt{\lambda_{\max}}}\Lambda^{\scriptscriptstyle 1/2}$. Then, we have
\begin{align*}
    &\zhadamarddopt(A,B,n-s) + s\log\lambda_{\max}\\
   &\quad = \ldet(B^\top B) + \textstyle\sum_{\ell=1}^{n-s} \log \delta_{\ell}(I + A(B^\top B)^{-1}A^\top) + s\log\lambda_{\max}\\
     &\quad = \ldet((1/\lambda_{\max})C) +\textstyle\sum_{\ell=1}^{n-s} \log \delta_{\ell}(\lambda_{\max} C^{-1}) + s\log\lambda_{\max}\\
     &\quad = \ldet(C) +\textstyle\sum_{\ell=1}^{n-s} \log \delta_{\ell}(C^{-1}). 
\end{align*}
\qed 
\end{proof}


\subsubsection{NLP-Id and Natural bounds.}\label{sec:NLP-IDandNatural}  In this subsection, we 
establish close connections between the NLP-Id bound for $\MESP$ and the
natural bound for 0/1 $\DOPT$.

\begin{lemma}\label{lem:natural_to_nlp_modified}
Let $C\in\mathbb{S}^n_{+}$\,, $0<s<n$, $0\leq \hat x\leq \mathbf{e}$, with $\mathbf{e}^\top \hat{x}=s$, and  $0<\gamma\leq 1/\lambda_{\max}$\,. Define $\Theta(\gamma):= I_n + \Diag(\mathbf{e}-\hat{x})^{\scriptscriptstyle 1/2}(\gamma C-I_n)\Diag(\mathbf{e}-\hat{x})^{\scriptscriptstyle 1/2}$ and let 
     \begin{equation*} 
         h(\gamma) := \ldet( \Theta(\gamma)) - (n-s)\log\gamma.
     \end{equation*}
Assume that $\Theta(\gamma) \succ 0$ for all $\gamma \in \left(0,1/\lambda_{\max}\right]$.  Then $h(\gamma)$ is {nonincreasing} in $\gamma$ on $\left(0,1/\lambda_{\max}\right]$.
Moreover, if $\hat{x} \notin \{0,1\}^n$, then $h(\gamma)$ is {strictly decreasing} in $\gamma$ on $\left(0,1/\lambda_{\max}\right)$.
\end{lemma}

\begin{proof}
We have that 
\begin{equation}\label{derh}
\begin{array}{ll}
    h'(\gamma) &= \gamma^{-1}   \Tr(\Theta(\gamma)^{-1}\gamma \Diag(\mathbf{e}-\hat{x})^{\scriptscriptstyle 1/2} C\Diag(\mathbf{e}-\hat{x})^{\scriptscriptstyle 1/2} ) - \gamma^{-1}(n-s)\\
    &=\gamma^{-1}\Tr(\Theta(\gamma)^{-1}(\Theta(\gamma) - \Diag(\hat{x}))) -  \gamma^{-1}\Tr(\Diag(\mathbf{e}-\hat{x})) \\
     &=\gamma^{-1}\Tr(I_n -\Theta(\gamma)^{-1}\Diag(\hat{x})) -  \gamma^{-1}\Tr(\Diag(\mathbf{e}-\hat{x})) \\
    &= \gamma^{-1}\Tr(\Diag(\hat{x}) - \Theta(\gamma)^{-1}\Diag(\hat{x}))\\
    &= \gamma^{-1}\Tr((I_n - \Theta(\gamma)^{-1})\Diag(\hat{x}))\\
     &= \gamma^{-1}\Tr(\Diag(\hat{x})^{\scriptscriptstyle 1/2}(I_n-\Theta(\gamma)^{-1})\Diag(\hat{x})^{\scriptscriptstyle 1/2})\\
     &= \gamma^{-1}\sum_{j \in N} \hat{x}_j  (I_n-\Theta(\gamma)^{-1})_{jj}\,.
\end{array}
\end{equation}
As $\gamma \in (0, 1/\lambda_{\max}]$\,, then $\gamma C \preceq I_n \Rightarrow \Theta(\gamma) \preceq I_n \Rightarrow \Theta(\gamma)^{-1} \succeq I_n\,,$ so we conclude that $h$ is nonincreasing in the interval $(0,1/\lambda_{\max}]$. Now consider that $\hat{x} \not\in \{0,1\}^n$ and  $ \gamma \in (0,1/\lambda_{\max})$. Let $\hat\jmath \in N$, such that $\hat{x}_{\hat\jmath} \in (0,1)$. Then, 
\[
\Theta(\gamma)_{\hat\jmath\hat\jmath} = \hat{x}_{\hat\jmath} + (1-\hat{x}_{\hat\jmath})\gamma C_{\hat{\jmath}\hat{\jmath}} < 1,
\]
and because $\Theta(\gamma) \circ \Theta(\gamma)^{-1} \succeq I_n$ (see, e.g., \cite[Thm. 7.7.9(c)]{HJBook}), then we have that
\[
\Theta(\gamma)_{\hat{\jmath}\hat{\jmath}}(\Theta(\gamma)^{-1})_{\hat{\jmath}\hat{\jmath}}\geq 1 \Rightarrow(\Theta(\gamma)^{-1})_{\hat{\jmath}\hat{\jmath}} > 1.
\]
The result follows.
\qed 
\end{proof}

\begin{lemma}\label{lem:rankA}
    Let $A \in \mathbb{R}^{n\times m}$, $B \in \mathbb{R}^{q\times m}$, with
{\scriptsize$\begin{pmatrix}A\\ B \end{pmatrix}$} being full column rank.  
 Let $C:=I_n - A(A^\top A +B^\top B)^{-1}A^\top$ and $\lambda_{\max}$ be the greatest eigenvalue of $C$. Then, $\lambda_{\max}\leq 1$, with equality if and only if $\rank(A)<n$. 
\end{lemma}
\begin{proof}
    Let $E := A^\top A +B^\top B$. We note that 
$\lambda_{\max} = 1 - \lambda_{n}(AE^{-1}A^\top)$  with $\lambda_{n}(AE^{-1}A^\top)\allowbreak\geq 0$,  then  
$\lambda_{\max}=1 \Leftrightarrow \lambda_{n}(AE^{-1}A^\top) = 0  \Leftrightarrow \rank(AE^{-1}A^\top) < n$. As $\rank(AE^{-1}A^\top) \allowbreak= \rank(AE^{\scriptscriptstyle -1/2}) = \rank(A)$, the result follows.
\qed \end{proof}

\begin{theorem}\label{thm:natural-nlp-equiv}
    Let $\MESP(C,n-s):= \mathcal{M}(\DOPT(A,B,s))$
    for a 0/1 $\DOPT$ instance $\DOPT(A,B,s)$. 
Then, the $\mathcal{M}$-induced  $\NLPId$ bound 
\begin{itemize}
    \item coincides with the natural bound for $\DOPT(A,B,s)$ if $\rank(A) < n$.
    \item strictly dominates the natural bound for $\DOPT(A,B,s)$ if  $\rank(A)=n$ and there exists a non-binary optimal solution to the $\mathcal{M}$-induced  $\NLPId$ bound.
\end{itemize}
\end{theorem}
\begin{proof}
By the definition of $\mathcal{M}$, we have 
$C := I_n - A(A^\top A +B^\top B)^{-1}A^\top$.  Let $\hat{x}$ be a feasible solution of \ref{natural_bound} with finite objective value. Then, we have 
\begin{align}
        &\ldet\left( A^\top \Diag(\hat{x})A + B^\top B \right)\label{relnatnlpid}\\
        &~= \ldet\left( A^\top A + B^\top B - A^\top\Diag(\mathbf{e}-\hat{x})A\right)\nonumber\\
        &~=\ldet\left( I_n - (A^\top A + B^\top B)^{\scriptscriptstyle -1/2}A^\top\Diag(\mathbf{e}-\hat{x})A(A^\top A + B^\top B)^{\scriptscriptstyle -1/2} \right) \nonumber\\[-2pt]
        &\qquad + \ldet(A^\top A + B^\top B)\nonumber\\
        &~= \ldet( I_n - \Diag(\mathbf{e}-\hat{x})^{\scriptscriptstyle 1/2}(A(A^\top A + B^\top B)^{-1}A^\top)\Diag(\mathbf{e}-\hat{x})^{\scriptscriptstyle 1/2} \
        )\nonumber\\[-2pt]
        &\qquad + \ldet(A^\top A + B^\top B)\nonumber\\
        &~= \ldet\left( I_n + \Diag(\mathbf{e}-\hat{x})^{\scriptscriptstyle 1/2}(C-I_n)\Diag(\mathbf{e}-\hat{x})^{\scriptscriptstyle 1/2} \right) \nonumber+  \ldet(A^\top A + B^\top B). 
\end{align}
If $\rank(A) < n$, then from Lemma \ref{lem:rankA}, we have  $\lambda_{\max} = 1$, and the result for this case follows.
 Consider the case  $\rank(A)=n$. From Lemma \ref{lem:rankA} we have $\lambda_{\max} < 1$. Let $x^*$ be an optimal solution for the natural bound for $\DOPT(A,B,s)$, and let $\mathbf{e}-x^*_{\gamma}$ be a non-binary optimal solution for the $\mathcal{M}$-induced  $\NLPId$ bound. Then,  
\begin{align*}
    &\ldet\left( I_n + \Diag(\mathbf{e}-x^*_{\gamma})^{\scriptscriptstyle 1/2}\big(\textstyle\frac{1}{\lambda_{\max}}C-I_n\big)\Diag(\mathbf{e}-x^*_{\gamma})^{\scriptscriptstyle 1/2} \right) +  \ldet(A^\top A + B^\top B)\\
    &~<\ldet\left( I_n + \Diag(\mathbf{e}-x^*_{\gamma})^{\scriptscriptstyle 1/2}(C-I_n)\Diag(\mathbf{e}-x^*_{\gamma})^{\scriptscriptstyle 1/2} \right) +  \ldet(A^\top A + B^\top B)\\
    &~\leq\ldet\left( I_n + \Diag(\mathbf{e}-x^*)^{\scriptscriptstyle 1/2}(C-I_n)\Diag(\mathbf{e}-x^*)^{\scriptscriptstyle 1/2} \right) +  \ldet(A^\top A + B^\top B),
\end{align*}
where the first inequality uses Lemma \ref{lem:natural_to_nlp_modified}  and the second uses optimality of $x^*$.
\qed \end{proof}

\begin{remark}\label{rem:pure-dopt-nlpid-natural}
For $A\in\mathbb{R}^{n\times m}$  with $\rank(A)=m<n$, consider the pure 0/1   $\DOPT$ instance $\DOPT(A,0,s)$.
    Let $\MESP(C,n-s):= \mathcal{M}(\DOPT(A,0,s))$.
   Then, from Theorem \ref{thm:natural-nlp-equiv}, we have that $\znatural(A,0,s)= \znlpidmesp(C,n-s) + \ldet(A^\top A).$
\end{remark}


In Theorem \ref{thm:nlp-natural-equiv}, we show that the bound for \ref{MESP} that is derived from the  natural bound for the $\DOPT$ problems established by $\mathcal{D}$  
is equivalent to the $\NLPId$ bound for \ref{MESP}. The following lemma is used in its proof.

\begin{lemma}\label{lem:nlp-natural-equiv-every-x}
    Let  $C \!\in \!\mathbb{S}^n_{+}$  with real Schur decomposition $\Phi\Lambda\Phi^\top$ and greatest eigenvalue $\lambda_{\max}$\,. Let $A:=\Phi(I_n-\frac{1}{\lambda_{\max}}\Lambda)^{\scriptscriptstyle 1/2}$ and $B:=\frac{1}{\sqrt{\lambda_{\max}}}\Lambda^{\scriptscriptstyle 1/2}$. Then, we have that
$
\ldet\!\left( I_n + \Diag(\hat{x})^{\scriptscriptstyle 1/2}((1/\lambda_{\max}) C\!-\!I_n)\Diag(\hat{x})^{\scriptscriptstyle 1/2} \right) = \ldet\!\left( {A}^\top \Diag(\mathbf{e}\!-\!\hat{x}) {A} \!+\! {B}^\top {B}\right),
$
for any $0\leq \hat{x}\leq \mathbf{e}$ such that  the $\ldet$ functions above have a finite value. 
\end{lemma}
\begin{proof}
   We have  
    \begin{align*}
        &\ldet\left( I_n + \Diag(\hat{x})^{\scriptscriptstyle 1/2}((1/\lambda_{\max}) C-I_n)\Diag(\hat{x})^{\scriptscriptstyle 1/2} \right) \\
        &\quad  = \ldet\left( I_n - \Diag(\hat{x})^{\scriptscriptstyle 1/2}{A}{A}^{\top}\Diag(\hat{x})^{\scriptscriptstyle 1/2}\right)
        =\ldet\left( I_n - {A}^\top\Diag(\hat{x}){A}\right)\\
        &\quad = \ldet\left( {A}^\top {A} +B^\top B   - {A}^\top \Diag(\hat{x}) {A}\right)
    =\ldet\left( {A}^\top \Diag(\mathbf{e}-\hat{x}) {A} + {B}^\top {B}\right),
    \end{align*}
where  we use that ${A}^\top {A} +B^\top B =I_n$\,.
\qed \end{proof}

\begin{theorem}\label{thm:nlp-natural-equiv}
    For  $C \!\in \!\mathbb{S}^n_{+}$\, with greatest eigenvalue $\lambda_{\max}$\,,  consider the $\MESP$ instance $\MESP(C,s)$ and the real Schur decomposition $\Phi\Lambda\Phi^\top$ of $C$. Let
$\DOPT(\allowbreak A, B, n-s):= \mathcal{D}(\MESP(C,s);\Phi)$. Then, 
the $\mathcal{D}$-induced natural bound is equal to the $\NLPId$  bound for $\MESP(C,s)$; that is, 
$\znatural(A,B,n-s) + s\log \lambda_{\max}=\znlpidmesp(C,s)$.
\end{theorem}

\begin{proof}  
   By the definition of $\mathcal{D}$, we have 
$A:=\Phi(I_n-\frac{1}{\lambda_{\max}}\Lambda)^{\scriptscriptstyle 1/2}$ and $B:=\frac{1}{\sqrt{\lambda_{\max}}}\Lambda^{\scriptscriptstyle 1/2}$.   Let $\hat{x}$ be a feasible solution to $\NLPId$ with finite objective value. Then, from Lemma \ref{lem:nlp-natural-equiv-every-x}, we have  
\begin{equation*} 
\begin{array}{l}
        \ldet\left( I_n + \Diag(\hat{x})^{\scriptscriptstyle 1/2}((1/\lambda_{\max}) C-I_n)\Diag(\hat{x})^{\scriptscriptstyle 1/2} \right) + s\log\lambda_{\max}\\ 
     \quad =\ldet\left( {A}^\top \Diag(\mathbf{e}-\hat{x}) {A} + {B}^\top {B}\right) +s\log\lambda_{\max}\,,
    \end{array}
    \end{equation*}
where on the left-hand side, we have the objective value of $\NLPId$ at $\hat{x}$,  and on the  right-hand side, we have the objective value of \ref{natural_bound} at the corresponding feasible solution $\mathbf{e}\!-\!\hat{x}$, added to  $s\!\log\!\lambda_{\max}$\,.
The result follows.
\qed \end{proof}


\subsubsection{NLP-Id and NLP-Tr.}

 As described in \S\ref{sec:MESPbounds}, there are three known strategies for selecting the parameters for the \ref{nlp_bound} bound for $\MESP$ instances, leading to the special cases $\NLPId$, $\NLPDi$ and $\NLPTr$. In Theorem \ref{thm:nlp_trace_eq_id}, we present sufficient conditions for $\NLPTr$ to be precisely $\NLPId$, when applied to the $\MESP$ instances that arise from
 \emph{pure} 0/1 $\DOPT$ instances, under the mapping $\mathcal{M}$. 

\begin{theorem}\label{thm:nlp_trace_eq_id}
    Let $C := I_n - UU^\top$, where $U^\top U = I_m$\,, with $m< n$, and consider the $\MESP$ instance $\MESP(C,s)$.  Suppose that there exists $\hat \beta \in \mathbb{R}^n_+$ such that $(C\circ C)\hat\beta = \mathbf{e}$. Then, we have that $\NLPTr$ is precisely  $\NLPId$. 
\end{theorem}
\begin{proof}
     The dual of \eqref{eq:nlp-tr}
is  
    \begin{equation}\label{eq:nlp-trace-dual}
        \max\,\{\Tr(C\Omega)\,:\, \diag(\Omega)=\mathbf{e},\, \Omega\succeq 0\}\,.
    \end{equation}
    Note that, because $C=I_n-UU^\top$ and $\diag(\Omega)=\mathbf{e}$ for every feasible solution $\Omega$ of  \eqref{eq:nlp-trace-dual},   we can rewrite the objective function  as $n - \Tr(U^\top \Omega U)$. Next, we will  present  a feasible  solution  $\widehat\Omega$ to \eqref{eq:nlp-trace-dual}, such that $U^\top \widehat\Omega U = 0$. 

    Let $\hat\beta$ be a solution of $(C\circ C)\beta = \mathbf{e}$ with $\hat\beta \geq \mathbf{0}$. Construct $\widehat\Omega := C\Diag(\hat\beta)C$.  We  note that $\widehat\Omega \succeq 0$, $\widehat\Omega_{ii} =  \textstyle\sum_{j = 1}^n  C_{i j}^2 \hat\beta_j = 1$ and $U^\top \widehat\Omega U = U^\top (I_n - UU^\top) \Diag(\hat\beta) (I_n-UU^\top) U = 0$. Thus, $\widehat\Omega$ is a feasible solution to \eqref{eq:nlp-trace-dual} with objective value equal to $n$.
    
    Now, we note that $\hat Y := I_n$ is a feasible solution to \eqref{eq:nlp-tr} with objective value equal to $n$ as well, and therefore, $ \hat Y:=I_n$ is an optimal solution to \eqref{eq:nlp-tr}. Next, we will show that  $\hat{Y}$ is the unique optimal solution. 
    
    Consider any optimal solution $Y$ to \eqref{eq:nlp-tr}, then $\Tr(Y) = n$. Note that 
    $C \widehat\Omega = \widehat\Omega$ because $C^2=C$. Then, we have 
    \[
    \Tr((Y-C)\widehat\Omega) = \Tr(Y\widehat\Omega) -  \Tr(\widehat\Omega) = \textstyle\sum_{i = 1}^n Y_{ii} \widehat\Omega_{ii} - n =   \textstyle\sum_{i = 1}^n Y_{ii} - n = 0.
    \]
    As $(Y\!-\!C) \succeq 0$, $\widehat\Omega \succeq 0$ and $\Tr((Y\!-C\!)\widehat\Omega) = 0$, we have $(Y\!-\!C)\widehat\Omega = 0$, which implies $(Y\!-\!I_n)\widehat\Omega = 0$, considering  that $C \widehat\Omega \!=\! \widehat\Omega$. Then, for all $i \!\in\! N$, we have $(Y_{ii}\!-\!1)\widehat\Omega_{ii}\!=\! 0$. As $\diag(\widehat\Omega) \!=\! \mathbf{e}$,  $\hat Y$ is the unique optimal solution to~\eqref{eq:nlp-tr}.

    Finally, we note that when using the $\Trace$ strategy to select the \ref{nlp_bound} parameters, if $D$ is selected as $I_n$\,, then we must select $\gamma=1$.
    Therefore, as $\lambda_{\max}=1$ for $C:=I_n-UU^\top$, both parameters $D$ and $\gamma$ are selected as in the $\Identity$ strategy, and we conclude that the $\Trace$ and $\Identity$ strategies lead to the same \ref{nlp_bound}.
\qed \end{proof}

With the examples below, we show that for all $n\geq 4$, the existence of a nonnegative solution to the linear system $(C\circ C)\beta=\mathbf{e}$ mentioned in the theorem above is sufficient, but not necessary, for $I_n$ to be an optimal solution to \eqref{eq:nlp-tr}.

\begin{example}
    Consider $u_1 :=  \frac{1}{\sqrt{110}}(7,6,3,4,0,0,\ldots,0)^\top$, $u_2 :=  \frac{1}{\sqrt{10}}(-1,2,1,-2,0,0,\allowbreak \cdots\!,0)^\top$, 
    $\omega := (1,-1,1,-1,1,1,\ldots,1)^\top$. Let $U := (u_1\; u_2)$, $\widehat\Omega := \omega\omega^\top$.
    Note that $\widehat\Omega$ is a feasible solution to \eqref{eq:nlp-trace-dual} with $U^\top \Omega U = (\omega^\top U)^\top(\omega^\top U)= 0$,  with objective value $n$. Therefore, $I_n$ is optimal for \eqref{eq:nlp-tr}. Consider $C:= I_n-UU^\top$ and $\Psi := C\circ C$. Let $y :=\textstyle\bigl(\frac{5}{8},-\frac{9}{2},\frac{19}{8},0,0,0,\ldots,0\bigr)^\top$. We can verify $\mathbf{e}^\top y = -\frac{3}{2}<0$ and $\Psi^\top y = \mathbf{e}_{3}\geq \mathbf{0}$.     Then, by the Farkas Lemma, there does not exist $\hat{\beta} \in \mathbb{R}^n_{+}$ with $\Psi \hat\beta = \mathbf{e}$.
\end{example}

Next, we will establish easily verifiable conditions under which the linear system $(C\circ C)\beta=\mathbf{e}$ has a nonnegative solution. We base our  proof 
on the Gauss-Seidel iterative method for linear systems and the (nonsingular) Sassenfeld matrices (see \cite{sass}, and more accessibly  \cite{baumann2017note} and \cite[Sec. 4.3]{Wendland_2017}).

\begin{definition}
    Given 
$\Psi\in\mathbb{R}^{n\times n}$, the \emph{Sassenfeld coefficient} $\alpha\in\mathbb{R}^n$ is defined by

    \vspace{-15pt}
    \begin{align}\label{sassenfeldcoef}
        \alpha_i := \begin{cases}
            \sum_{j=2}^n |\Psi_{1j}|/|\Psi_{11}|\,,\quad &i=1;\\[5pt]
           \left (\textstyle\sum_{j =1}^{i-1}\alpha_j |\Psi_{ij}| + \sum_{j=i+1}^n |\Psi_{ij}|\right)/|\Psi_{ii}|\,,\quad &i = 2,\dots,n.
        \end{cases}
    \end{align}
If $\alpha_i<1$, for all $i\in N$, then we say that $\Psi$ is a \emph{Sassenfeld matrix.}
\end{definition}

\begin{theorem}\label{thm:sassenfeld}
     Let $\Psi\in\mathbb{R}^{n\times n}$ be a nonnegative Sassenfeld matrix, such that\break $\sum_{\substack{j\in N\\j\neq i}}\frac{\Psi_{ij}}{\Psi_{jj}} \leq 1$, for all $i\in N$. 
Then, 
the linear system $\Psi \beta = \mathbf{e}$, has a unique solution $\hat\beta$, with $\hat\beta\geq \mathbf{0}$. 
\end{theorem}

\begin{proof}
    As $\Psi$ is a Sassenfeld matrix, it is nonsingular, and for any given $\beta^0\in\mathbb{R}^n$, the iterates of the Gauss-Seidel method, defined by 
    \begin{equation}\label{gauss-seidel}
    \textstyle \beta^{k+1}_i := \frac{1}{\Psi_{ii}}\left(1-\sum_{\substack{j\in N\\j< i}} \Psi_{ij}\beta_j^{k+1}
    - \sum_{\substack{j\in N\\j> i}} \Psi_{ij}\beta_j^k  
    \right),\quad i \in N,
    \end{equation}
    for $k \in \mathbb{N}$, 
converge to the unique solution of the linear system $\Psi\beta=\mathbf{e}$ (see \cite[Thms. 2 and 3]{baumann2017note}). Next, we will prove that the unique solution is nonnegative, under the additional conditions $\Psi\geq 0$ and  $\sum_{\substack{j\in N\\j\neq i}}\frac{\Psi_{ij}}{\Psi_{jj}} \leq 1$, for all $i\in N$.

  We will prove by induction (on $k$) that if $\beta^0 := \mathbf{0}$, then for any $k \in \mathbb{N}$, we have $0 \leq \beta_i^k \leq 1/\Psi_{ii}$\,, for all $i \in N$. The upper bound on $\beta^k$ will be used to assure that all iterates are  nonnegative. 
    For the induction, the base case ($k=0$)  is trivial. Now assume  that $0 \leq \beta_i^k\leq
    1/\Psi_{ii}$\,, for all $i\in N$. Then, we should prove that $0 \leq \beta_i^{k+1} \leq
    1/\Psi_{ii}$\,, for all $i\in N$; we will prove it by induction on $i$, using the nonnegativity of $\Psi$.
    For $i=1$, using the upper bounds on $\beta^k$, we see that  
     \begin{equation*}
        0\leq \textstyle   \sum_{\substack{j\in N\\j>1}}\Psi_{1j}\beta_j^{k} \leq \textstyle \sum_{\substack{j\in N\\j\neq 1}}\frac{\Psi_{1j}}{\Psi_{jj}} \leq 1.
    \end{equation*}    
    Then, from \eqref{gauss-seidel}, we see that  $0\leq \beta^{k+1}_1\leq 1/\Psi_{11}$\,.
    Now, assume that, for $i=2,\ldots, \ell-1$, we have
\begin{align}\label{inductioni}
        0&\leq \textstyle   \sum_{\substack{j\in N\\j< i}} \Psi_{ij}\beta_j^{k+1} + \sum_{\substack{j\in N\\j>i}}\Psi_{ij}\beta_j^{k} \leq 1.
    \end{align}
     Then, from \eqref{gauss-seidel}, we see that   $0\leq \beta^{k+1}_i\leq 1/\Psi_{ii}$\,, for $i=2,\ldots, \ell -1$.
     
    Next  we prove \eqref{inductioni} for $i=\ell$, using  bounds on $\beta^k$ and  $\beta^{k+1}_i$ ($i<\ell$). We have
    \begin{equation*}
        0\leq \textstyle\sum_{\substack{j\in N\\j<\ell}}\Psi_{\ell j}\beta_j^{k+1} +   \sum_{\substack{j\in N\\j>\ell}}\Psi_{\ell j}\beta_j^k\leq  \sum_{\substack{j\in N\\j\neq \ell}}\frac{\Psi_{\ell j}}{\Psi_{jj}}\leq 1.
    \end{equation*}
Therefore, $0\leq \beta^{k+1}_\ell\leq 1/\Psi_{\ell \ell}$\,.

From  induction on $i$, we have $0\leq \beta^{k+1}_i\leq 1/\Psi_{ii}$\,, for all $i\in N$.
By induction on $k$, we conclude that $\beta^k\geq \mathbf{0}$, for all $k\in\mathbb{N}$, and because $\{\beta \in \mathbb{R}^n_+ ~:~ \beta_\ell \leq 1/\Psi_{\ell\ell} \mbox{ for } \ell \in N\}$  
is a compact set,
we can finally conclude that the Gauss-Seidel method converges to a nonnegative solution of the linear system $\Psi\beta=\mathbf{e}$.
\qed \end{proof}
Next we see 
the assumption  
$\textstyle\sum_{\substack{j\in N\\j\neq i}}\frac{\Psi_{ij}}{\Psi_{jj}} \leq 1$, for $i\!\in \!N$, in Theorem \ref{thm:sassenfeld} is necessary. 

\begin{example}
    Let $C:=I_n-UU^\top$, where $U := \frac{1}{\sqrt{14}}(1,2,3,0,0,\ldots,0)^\top $. Let $\Psi:=C\circ C$ and consider 
    $\alpha_i$, for $i\in N$, defined by \eqref{sassenfeldcoef}.   We have that  $\max_{i \in N}\alpha_i \leq 0.56$, $\max_{i \in N}\sum_{\substack{j\in N\\j\neq i}}\frac{\Psi_{ij}}{\Psi_{jj}} \approx 1.46$,  
    $\Psi$ is full-rank, and  
    $\hat\beta :=(2/3,-5/3,10,1,1,\ldots,1)^\top$     
    is the unique solution of $\Psi\beta = \mathbf{e}$.
\end{example}

Next, we establish a stronger and even more-easily-verifiable condition under which we guarantee that the linear system $\Psi\beta=\mathbf{e}$ has a nonnegative solution, where $\Psi:=C\circ C$, $C:=I_n-UU^\top$, and $U^\top U=I_m$\,, with $m<n$.

\begin{lemma}\label{lem:CoC_dd_strict}
    Let $C := I_n - UU^\top$, where $U^\top U = I_m$\,, with $m< n$. Let $\Psi:=C\circ C$. If $C_{ii} > 1/2$ for all $i \in N$, then  $\Psi$ is  strictly diagonally dominant. 
\end{lemma}
\begin{proof}
    Note that 
    $C^2= C$.
    So, for all $i\in N$, 
    we have 
    $ C_{ii} = C_{i\cdot} C_{i\cdot}^{\top} = \sum_{j\in N} C_{ij}^2\,.$
     Then, for all $i  \in N$, we have  
    \begin{equation}\label{eq:sum-off-diag-psi-dd}
        \textstyle \sum_{\substack{j\in N\\j\neq i}} |\Psi_{ij}|  = \textstyle\sum_{\substack{j\in N\\j\neq i}} \Psi_{ij} =  \sum_{j\in N} C_{ij}^2 -C_{ii}^2  = C_{ii} -C_{ii}^2\,.
    \end{equation}
   Therefore,
    $
    \textstyle |\Psi_{ii}| - \sum_{\substack{j\in N\\j\neq i}} |\Psi_{ij}| = 2C_{ii}^2 - C_{ii}\,.
    $
    The result follows.
\qed \end{proof}

\begin{proposition}\label{prop:diagC}
     Let $C := I_n - UU^\top$, where $U^\top U = I_m$\,, with $m<n$. Let $\Psi:=C\circ C$. If $C_{ii} > 1/2$ for all $i \in N$,  then  $\Psi$ is a  Sassenfeld matrix and  $\sum_{\substack{j\in N\\j\neq i}}\frac{\Psi_{ij}}{\Psi_{jj}} \leq 1$, for all $i\in N$. 
\end{proposition}

\begin{proof}
From Lemma \ref{lem:CoC_dd_strict}, we have that $\Psi$ is strictly diagonally dominant.
Then, $\sum_{\substack{j\in N\\j\neq i}} |\Psi_{ij}|/|\Psi_{ii}|<1$, for all $i\in N$, and by induction on $i$ we can easily verify that $\Psi$ is a Sassenfeld matrix.
Furthermore, considering \eqref{eq:sum-off-diag-psi-dd} and $C_{ii}>1/2$ (and hence $\Psi_{ii}>1/4$), for all $i\in N$, we have that 
$\textstyle \sum_{\substack{j\in N\\j\neq i}}\frac{\Psi_{ij}}{\Psi_{jj}} \leq 4(C_{ii} - C_{ii}^2)\leq 1$.   
\qed \end{proof}

Finally, the family of examples below shows that, for all $n\geq 4$, 
the conditions stated in Theorem \ref{thm:sassenfeld},
 where $\Psi:=C\circ C$,  $C:=I_n-UU^\top$, and $U^\top U=I_m$, with $m<n$, 
  do not imply that $C_{ii}>1/2$ for all $i\in N$.

\begin{example}
    Let $C:=I_n-UU^\top$, where $U := \frac{1}{\sqrt{30}}(1,2,3,4,0,0,\ldots,0)^\top $. Let $\Psi:=C\circ C$ and consider the Sassenfeld coefficients $\alpha_i$, for $i\in N$, defined by \eqref{sassenfeldcoef}.  We have $\max_{i \in N}\alpha_i \leq 0.34$, $\max_{i \in N}\sum_{\substack{j\in N\\j\neq i}}\frac{\Psi_{ij}}{\Psi_{jj}} \leq 0.8$, and $\diag(C)=\textstyle\bigl(\frac{29}{30},\frac{26}{30},\frac{21}{30},\frac{14}{30},1,1,\ldots\allowbreak ,1\bigr)^\top$.
 Illustrating again the result of Theorem \ref{thm:sassenfeld},  $\Psi$ is full rank, and  the nonnegative vector 
    $\hat\beta :=\textstyle\bigl(\frac{918}{929},\frac{873}{929},\frac{698}{929},\frac{3393}{929},1,1,\ldots,1\bigr)^\top$     
    is the unique solution of $\Psi \beta = \mathbf{e}$.
\end{example}

\begin{corollary}\label{cor:jacobi}

Let $C := I_n - UU^\top$, where $U^\top U = I_m$\,, with $m< n$, and consider the $\MESP$ instance $\MESP(C,s)$.  Let $\Psi:=C\circ C$ and assume that any of the following successively weaker conditions holds:
    \begin{itemize}
        \item $C_{ii} > 1/2$, for all $i \in N$;\\[-10pt]
        \item $\Psi$ is a Sassenfeld matrix with $\sum_{\substack{j\in N\\j\neq i}}\frac{\Psi_{ij}}{\Psi_{jj}} \leq 1$, for all $i\in N$;\\[-6pt]
        \item $\Psi\beta=\mathbf{e}$ has a nonnegative solution.
    \end{itemize}
        Then $D=I_n$ is the unique solution of the trace minimization problem in \eqref{eq:nlp-tr}, and  we have that $\NLPTr$ is precisely  $\NLPId$.
\end{corollary}

\begin{proof}
The result follows from Theorems \ref{thm:nlp_trace_eq_id} and
\ref{thm:sassenfeld}, and Proposition \ref{prop:diagC}.   
\qed \end{proof}


\section{B\&B for MESP and 0/1 D-optimality with transferred upper bounds}\label{sec:bb}

In \S\ref{sec:compare-bounds}, we discussed the possibility of transferring upper bounds from instances of $\MESP$ to instances of  0/1 $\DOPT$ and vice-versa. Here, we address possible ways to apply this idea within a B\&B algorithm.

In Table \ref{tab:fixingBB}, we show the defining parameters of subproblems of interest in our following analysis, generated from the  instances of $\MESP$ and 0/1 $\DOPT$ when indices of $N$ are fixed in and out of the solution. 

\renewcommand{\arraystretch}{1.75} 
\begin{table}[!ht]
\centering
\begin{tabular}{l|l|l}
              & \multicolumn{1}{c|}{$\MESP(C,n-s)$}                   & \multicolumn{1}{c}{$\DOPT(A,B,s)$}        \\ \hline
Fix $j\in N$ out & $\MESP(C_{N\setminus  j ,N\setminus  j },n-s)$                           & $\DOPT(A_{(N\setminus  j )\cdot},B,s)$                \\ \hline
Fix $i\in N$ in  & $\MESP(C / C_{ii},n-s-1) + \log(C_{ii})$ & $\DOPT(A_{(N \setminus i )\cdot},\left({\begin{smallmatrix}\!\!A_{i\cdot}\!\!\\ B \end{smallmatrix}}\right),s-1)$
\end{tabular}
\caption{B\&B subproblems}\label{tab:fixingBB}
\end{table}


\subsection{Aiming to transfer bounds from MESP to 0/1 D-Opt}\label{subsec:M2D}

Considering 
Theorem \ref{thm:equivM}, we see that given an arbitrary 0/1 $\DOPT$ instance $\DOPT(A,B,s)$, 
we have that 
$\DOPT(A,B,s)\circeq\mathcal{M}(\DOPT(A,B,s))+\ldet(A^\top A + B^\top B)$, so $\zdoptthing(A,B,s)= \zmespthing(\mathcal{M}(\DOPT(A,B,s)))+\ldet(A^\top A + B^\top B)$. 
 
Therefore, when applying B\&B to $\DOPT(A,B,s)$, aiming to use upper bounds developed for $\MESP$, we can apply $\mathcal{M}(\cdot)$ to each  0/1 $\DOPT$ subproblem generated, and then compute an upper bound for each $\MESP$ problem obtained, with its objective function added to the constant $\ldet(\tilde{A}^\top\tilde{A} + \tilde{B}^\top\tilde{B})$, where $\tilde{A}$ and $\tilde{B}$ are the input matrices for the  0/1 $\DOPT$ subproblem.

Another way to solve $\DOPT(A,B,s)$ using upper bounds developed for $\MESP$, is to apply B\&B directly to $\mathcal{M}(\DOPT(A,B,s))+\ldet(A^\top A + B^\top B)$ (this procedure was applied in \cite{PonteFampaLeeMPB}).

In Theorem \ref{thm:bbdopt} we show that when applying  B\&B  to the 0/1 $\DOPT$ instance $\DOPT(A,B,s)$ and applying  $\mathcal{M}(\cdot)$ to each  0/1 $\DOPT$ subproblem generated, the $\MESP$ subproblems created during the B\&B execution are exactly the same as when applying B\&B directly to $\mathcal{M}(\DOPT(A,B,s))+\ldet(A^\top A + B^\top B)$.
More specifically, we show that if we fix $i\in N$ in the solution ($j\in N$ out of the solution) of $\DOPT(A,B,s)$, then  the $\MESP$ problem obtained via $\mathcal{M}$, from the resulting 0/1 $\DOPT$ subproblem, is equivalent to the $\MESP$ subproblem obtained from $\mathcal{M}(\DOPT(A,B,s))+\ldet(A^\top A + B^\top B)$ by fixing $i\in N$ out of its solution ($j\in N$ in its solution).

For a given   0/1 $\DOPT$ instance $\DOPT(A,B,s)$ and $\MESP(C,n-s):= \mathcal{M}(\DOPT(\allowbreak A,B,s))$,  we consider the  following four problems, for any $i,j\in N$. 
\begin{equation}\label{fixinD}\tag{$i$-in-$\DOPT(A,B,s)$}
            \MESP(C_1,n-s) + \ldet\left(A_{(N \setminus i )\cdot}^\top A_{(N \setminus  i )\cdot} + \left({\begin{smallmatrix}\!\!A_{i\cdot}\!\!\\ B \end{smallmatrix}}\right)^\top \left({\begin{smallmatrix}\!\!A_{i\cdot}\!\!\\ B \end{smallmatrix}}\right)\right),
            \end{equation}
            where $ \MESP(C_1,n-s) := \mathcal{M}(\DOPT(A_{(N \setminus i )\cdot},\left({\begin{smallmatrix}\!\!A_{i\cdot}\!\!\\ B \end{smallmatrix}}\right),s-1))$. The subproblem \ref{fixinD}  is considered in a B\&B enumeration tree to solve $\DOPT(A,B,s)$, when we fix $i\in N$ in the solution on $\DOPT(A,B,s)$, and then create a $\MESP$ instance from the resulting 0/1 $\DOPT$ subproblem via $\mathcal{M}$.
            \begin{equation}\label{fixoutM} \tag{\mbox{$i$-out-$\mathcal{M}(\DOPT(A,B,s))$}}
            \MESP(C_{N\setminus i ,N\setminus  i },n-s) + \ldet(A^\top A + B^\top B).
            \end{equation}
            The subproblem \ref{fixoutM}  is considered in a B\&B enumeration tree to solve $\DOPT(A,B,s)$, when we  create a $\MESP$ instance from the given  0/1 $\DOPT$ instance  $\DOPT(A,B,s)$ via $\mathcal{M}$ and then fix $i\in N$ out of the solution  of the resulting $\MESP$ problem.
            \begin{equation}\label{fixoutD}\tag{$j$-out-$\DOPT(A,B,s)$}
            \MESP(C_0,n\!-\!1\!-\!s) + \ldet(A_{(N\setminus  j )\cdot}^\top A_{(N\setminus  j )\cdot} + B^\top B),
            \end{equation}
            where $\MESP(C_0,n-1-s):= \mathcal{M}(\DOPT(A_{(N\setminus  j )\cdot}~,B,s))$. The subproblem\break \ref{fixoutD}  is considered in a B\&B enumeration tree to solve $\DOPT(\allowbreak A,B,s)$, when we fix $j\in N$ out of the solution  on $\DOPT(A,B,s)$, and then create a $\MESP$ instance from the resulting 0/1 $\DOPT$ subproblem via $\mathcal{M}$.
            \begin{equation}\label{fixinM}\tag{\mbox{$j$-in-$\mathcal{M}(\DOPT(A,B,s))$}}
            \MESP(C / C_{jj},n\!\!-\!\!1\!-\!\!s)  + \log(C_{jj}) + \ldet(A^\top A + B^\top B).
            \end{equation}
           The subproblem \ref{fixinM}  is considered in a B\&B enumeration tree to solve $\DOPT(A,B,s)$, when we  create a $\MESP$ instance from the given  0/1 $\DOPT$ instance  $\DOPT(A,B,s)$ via $\mathcal{M}$ and then fix $j\in N$ in the solution  of the resulting $\MESP$ problem.

\begin{theorem}\label{thm:bbdopt}
Let $\DOPT(A,B,s)$ be  an arbitrary 0/1 $\DOPT$ instance and $\MESP(\allowbreak C,n-s):= \mathcal{M}(\DOPT(A,B,s))$. For any $i,j\in N$, we have that 
\begin{enumerate}
\item   \label{inDoutM}  \eqref{fixinD} and \eqref{fixoutM} are  exactly the same instance of $\MESP$ with the same constant added to their objective functions.
\item  \label{outDinM}   \eqref{fixoutD} and  \eqref{fixinM} are  exactly the same instance of $\MESP$ with the same constant added to their objective functions.
\end{enumerate}

\end{theorem}

\begin{proof}
       By the definition of $\mathcal{M}$,  we have that 
       \begin{align*}
           &C:=I_n - A(A^\top A+ B^\top B)^{-1} A^\top,\\ 
            &C_1 := I_{n-1} - A_{(N \setminus i )\cdot}(A_{(N \setminus i )\cdot}^\top ~A_{(N \setminus i )\cdot} +
\left({\begin{smallmatrix}\!\!A_{i\cdot}\!\!\\ B \end{smallmatrix}}\right)^\top \left({\begin{smallmatrix}\!\!A_{i\cdot}\!\!\\ B \end{smallmatrix}}\right))^{-1}A_{(N \setminus i )\cdot}^\top~,\\
            &C_0 := I_{n-1} - A_{(N\setminus  j )\cdot}(A_{(N\setminus  j )\cdot}^\top A_{(N\setminus  j )\cdot} +B^\top B)^{-1}A_{(N\setminus  j )\cdot}^\top~.
       \end{align*}
       
To prove Item  \ref{inDoutM}, we consider that $A_{(N \setminus i )\cdot}^\top~ A_{(N \setminus  i )\cdot} + \left({\begin{smallmatrix}\!\!A_{i\cdot}\!\!\\ B \end{smallmatrix}}\right)^\top \left({\begin{smallmatrix}\!\!A_{i\cdot}\!\!\\ B \end{smallmatrix}}\right) = A^\top A + B^\top B$ to 
 verify that $C_1\!=\!C_{N\setminus i ,N\setminus  i }$ and $\ldet\left(A_{(N \setminus i )\cdot}^\top A_{(N \setminus  i )\cdot} + \left({\begin{smallmatrix}\!\!A_{i\cdot}\!\!\\ B \end{smallmatrix}}\right)^\top \left({\begin{smallmatrix}\!\!A_{i\cdot}\!\!\\ B \end{smallmatrix}}\right)\right)=\ldet(A^\top A~+~ B^\top B)$. 

To prove Item \ref{outDinM}, we should verify that $C_0=C/C_{jj}$ and $\ldet(A_{(N\setminus  j )\cdot}^\top A_{(N\setminus  j )\cdot} + B^\top B)= \log(C_{jj})+\ldet(A^\top A + B^\top B)$. 
Let 
\[
        H\!:=\! \begin{pmatrix}\underset{\scriptscriptstyle (m+1)\times (m+1)}{X} & \underset{\scriptscriptstyle (m+1)\times (n-1)}{Y}\\[1.5ex]
\underset{\scriptscriptstyle (n-1)\times (m+1)}{Z} & \underset{\scriptscriptstyle (n-1)\times (n-1)}{W}
\end{pmatrix} 
\!:=\!
\begin{pmatrix}
    \begin{matrix}
       \underset{\scriptscriptstyle m\times m}{X_{11}} & \underset{\scriptscriptstyle m\times 1}{X_{12}}\\
        \underset{\scriptscriptstyle 1\times m}{X_{21}} & \underset{\scriptscriptstyle 1\times 1}{X_{22}}
    \end{matrix} & Y\\
    Z & W
\end{pmatrix} 
\!:=\!  
\begin{pmatrix}
    A^\top A+B^\top B &~~A_{j\cdot}^\top~~ & ~~A_{(N\setminus j )\cdot}^\top~~ \\
    A_{j\cdot} & 1 & \mathbf{0}^\top\\
    A_{(N\setminus  j )\cdot}& \mathbf{0} & I_{n-1}
\end{pmatrix}.
\]
 We have  
 $ H/X = I_{n-1} - 
      A_{(N\setminus  j )\cdot} (X^{-1})_{11} 
      A_{(N\setminus  j )\cdot}^\top$\,, and  (see, e.g., \cite[p.20]{SchurBook})
   \[
   (X^{-1})_{M,M} = (X/1)^{-1} = (A^\top A + B^\top B - A_{j\cdot}^\top A_{j\cdot})^{-1} = (A_{(N\setminus  j )\cdot}^\top  A_{(N\setminus  j )\cdot}  + B^\top B)^{-1},
   \]
   where $M$ is the ordered set 
   $\{1,\ldots,m\}$. 
   So, we conclude that $H/X = C_0$\,.

Also, we note that  
 $H/X_{11} = I_n - 
 \left({\begin{smallmatrix}A_{j\cdot}\\
A_{(N\setminus  j )\cdot} \end{smallmatrix}}\right)
 (A^\top A+ B^\top B)^{-1}(A_{j\cdot}^\top ~ A_{(N\setminus  j )\cdot}^\top) = \left({\begin{smallmatrix}C_{jj} & C_{j,N\setminus  j }\\
C_{N\setminus  j ,j} & C_{N\setminus  j ,N\setminus  j }\end{smallmatrix}}\right)$,
and $X/X_{11} = 1 - A_{j\cdot}(A^\top A + B^\top B)^{-1}A_{j\cdot}^\top = C_{jj}$\,.  
Then, as $H/X = (H/X_{11})/(X/X_{11})$ (see, e.g., \cite[Thm. 1.4]{SchurBook}),  we have $C_0 = C/C_{jj}$\,.

Additionally, we have that
  \begin{align*}
     & \ldet(A^\top A + B^\top B) + \log(C_{jj}) \\
     &\quad= \ldet(A^\top A + B^\top B) + \log(1 - A_{j\cdot}(A^\top A + B^\top B)^{-1}A_{j\cdot}^\top)\\
      &\quad= \ldet(A^\top A + B^\top B) + \ldet(I_m- (A^\top A + B^\top B)^{\scriptscriptstyle -1/2}A_{j\cdot}^\top A_{j\cdot}(A^\top A + B^\top B)^{\scriptscriptstyle -1/2})\\
      &\quad= \ldet(A^\top A + B^\top B-A_{j\cdot}^\top A_{j\cdot})
      \quad= \ldet(A_{(N\setminus  j )\cdot}^\top A_{(N\setminus  j )\cdot} + B^\top B).
  \end{align*}
\qed \end{proof}


\subsection{Aiming to transfer bounds from 0/1 D-Opt to MESP}

From 
Theorem \ref{thm:equivD}, given an arbitrary $\MESP$ instance $\MESP(\allowbreak C,s)$ and the real Schur decomposition $\Phi\Lambda\Phi^\top$ of $C$, 
we have  
$\MESP(C,s)\circeq\mathcal{D}(\MESP(C,s);\Phi)+s\log \lambda_{\max}$\,, so $\zmespthing(C,s)= \zdoptthing(\mathcal{D}(\MESP(C,s);\Phi))  +  s\log \lambda_{\max}$\,. 
Therefore, when applying B\&B to $\MESP(C,s)$, aiming to use upper bounds developed for 0/1 $\DOPT$, we can apply $\mathcal{D}(\cdot)$ to each  $\MESP$ subproblem, and then compute an upper bound for each 0/1 $\DOPT$ problem, with its objective function added to the constant $\tilde{s}\log \lambda_{1}(\tilde{C})$, where $\tilde{C}$ and $\tilde{s}$ are the input data  for the   $\MESP$ subproblem.

Similarly to what was explained in \S\ref{subsec:M2D}, there is another way to solve $\MESP(C,s)$ using upper bounds developed for 0/1 $\DOPT$, by  applying B\&B directly to $\mathcal{D}(\allowbreak \MESP(\allowbreak C,s);\Phi)+s\log \lambda_{\max}$\,.
Nevertheless,  differently from the analysis in \S\ref{subsec:M2D}, we will show in the following that when applying  B\&B  to the $\MESP$ instance $\MESP(C,s)$ and applying  $\mathcal{D}(\cdot)$ to each  $\MESP$ subproblem generated, the 0/1 $\DOPT$ subproblems created during the B\&B execution are \emph{not}  the same as when applying B\&B directly to $\mathcal{D}(\MESP(C,s);\Phi)+s\log \lambda_{\max}$\,.
More specifically, we show that if we fix $i\in N$ in the solution ($j\in N$ out of the solution) of $\MESP(C,s)$, then  the 0/1 $\DOPT$ problem obtained via $\mathcal{D}$, from the resulting  $\MESP$ subproblem, is \emph{not} equivalent to the 0/1 $\DOPT$ subproblem obtained from $\mathcal{D}(\MESP(C,s);\Phi)+s\log \lambda_{\max}$ by fixing $i\in N$ out of its solution ($j\in N$ in its solution).


Considering a given  $\MESP$ instance $\MESP(C,s)$,  the real Schur decomposition $\Phi\Lambda\Phi^\top$ of $C$,  and $\DOPT(A,B,n-s):= \mathcal{D}(\MESP(C,s);\Phi)$,  we are interested in the  following four problems, for any $i,j\in N$. 
\begin{equation}\label{fixinMi}\tag{\mbox{$i$-in-$
\MESP(C,s)$}}
            \DOPT(A_1,B_1,n-s) + (s\!-\!1)\log\lambda_{1}(C/C_{ii})  + \log(C_{ii}),
            \end{equation}
            where $\DOPT(A_1,B_1,n\!-\!s) := \mathcal{D}(\MESP(C/C_{ii},s\!-\!1);\tilde{\Phi})$ and  $\tilde{\Phi}\tilde{\Lambda}\tilde{\Phi}^\top$ is a  real~Schur decomposition of $C/C_{ii}$\,. 
The subproblem \ref{fixinMi}  is considered in a B\&B enumeration tree for $\MESP(C,s)$, when we fix $i\in N$ in the solution on $\MESP(C,s)$, and then create a 0/1 $\DOPT$ instance from the  $\MESP$  subproblem via $\mathcal{D}$.
            \begin{equation}\label{fixoutDi} \tag{\mbox{$i$-out-$\mathcal{D}(\MESP(C,s))$}}
            \DOPT(A_{(N\setminus  i )\cdot},B,n-s) + s\log\lambda_{\max}\,.
            \end{equation}
            The subproblem \ref{fixoutDi}  is considered in a B\&B enumeration tree for $\MESP(C,s)$, when we  create a 0/1 $\DOPT$ instance from the   $\MESP$ instance  $\MESP(C,s)$ via $\mathcal{D}$ and then fix $i\in N$ out of the solution  of the  0/1 $\DOPT$ problem.
            \begin{equation}\label{fixoutMj}\tag{\mbox{$j$-out-$
\MESP(C,s)$}}
            \DOPT(A_0,B_0,n\!-\!1\!-\!s) + s\log\lambda_{1}(C_{N\setminus  j ,N\setminus  j }),
            \end{equation}
            where $\DOPT(A_0,B_0,n\!-\!1\!-\!s):= \mathcal{D}(\MESP(C_{N\setminus  j ,N\setminus  j },s);\tilde{\Phi})$ and  $\tilde{\Phi}\tilde{\Lambda}\tilde{\Phi}^\top$ is a  real Schur decomposition of $C_{N\setminus  j ,N\setminus  j }$\,. 
            The subproblem \ref{fixoutMj}  is considered in a B\&B enumeration tree to solve $\MESP(C,s)$, when we fix $j\in N$ out of the solution  on $\MESP(C,s)$, and then create a 0/1 $\DOPT$ instance from the resulting $\MESP$ subproblem via $\mathcal{D}$.
            \begin{equation}\label{fixinDj}\tag{\mbox{$j$-in-$
\mathcal{D}(\MESP(C,s))$}}
           \DOPT(A_{(N \setminus j )\cdot},\left({\begin{smallmatrix}\!\!A_{j\cdot}\!\!\\ B \end{smallmatrix}}\right),n-s-1) + s \log\lambda_{\max}\,.
            \end{equation}
           The subproblem \ref{fixinDj}  is considered in a B\&B enumeration tree for $\MESP(C,s)$, when we  create a 0/1 $\DOPT$ instance from the   $\MESP$ instance  $\MESP(C,s)$ via $\mathcal{D}$ and then fix $j\in N$ in the solution  of the  0/1 $\DOPT$ problem.

With the following example, we show that, unlike what we proved in \S\ref{subsec:M2D} for the case where we transform 0/1 $\DOPT$ problems into $\MESP$ problems within a B\&B algorithm to solve 0/1 $\DOPT$, when we transform $\MESP$ problems into 0/1 $\DOPT$ problems within a B\&B algorithm to solve $\MESP$, we obtain different subproblems with the two strategies described above. In fact, we show not only that \ref{fixinMi} (resp., \ref{fixoutMj}) is not the same 0/1 $\DOPT$ instance as \ref{fixoutDi} (resp., \ref{fixinDj}), but also a stronger result. We show that the natural bound for these two instances are different,  the one calculated for the first being better (smaller) than the one calculated for the second.

\begin{example}
   Consider the $\MESP$ instance $\MESP(C,s)$, where $C:=\left({\begin{smallmatrix}3 & 2 & 0 \\
       2 & 2 & 0\\
       0 & 0 & 1 \end{smallmatrix}}\right)$.
    \begin{enumerate}
        \item Consider $s := 2$, $i :=1$. Then,
        \begin{itemize}
             \item $\znatural(A_1,B_1,n-s)  + (s-1)\log\lambda_{1}(C/C_{ii}) +  \log(C_{ii})\approx 1.099$.
    \item $\znatural(A_{(N\setminus  i )\cdot},B,n-s) + s\log\lambda_{\max} \approx 1.570$.
        \end{itemize}
        \item Consider $s := 1$, $j :=1$. Then, 
        \begin{itemize}
            \item $\znatural(A_0,B_0,n\!-\!1\!-\!s) + s\log\lambda_{1}(C_{N\setminus  j ,N\setminus  j }) \approx 0.693$.
            \item $\znatural\left(A_{(N \setminus j )\cdot},\left({\begin{smallmatrix}\!\!A_{j\cdot}\!\!\\ B \end{smallmatrix}}\right),n-s-1\right) + s \log\lambda_{\max} \approx 0.754$.
        \end{itemize}
    \end{enumerate}
\end{example}

We next establish that the behavior illustrated in this example always holds; i.e., the natural bound for \ref{fixinMi} (resp., \ref{fixoutMj}) is never greater than the natural bound for \ref{fixoutDi} (resp.,\break \ref{fixinDj}). This is useful, as it shows that when  transforming $\MESP$ problems into 0/1 $\DOPT$ problems within a B\&B algorithm to solve $\MESP$, with the goal of computing the natural bound for the 0/1 $\DOPT$ subproblems, it is better to branch on the $\MESP$ problem and transform each $\MESP$ subproblem into a 0/1 $\DOPT$ problem via $\mathcal{D}$ than transforming the given $\MESP$ problem into a 0/1 $\DOPT$ problem by applying $\mathcal{D}$ only once and then branching on the 0/1 $\DOPT$ problem. This observation is particularly relevant for those intending to apply B\&B based on the $\NLPId$ bound to $\MESP$. Although we have proved that the $\NLPId$ bound for a $\MESP$ instance is the same as the natural bound for the 0/1 $\DOPT$ instance derived from it via $\mathcal{D}$, it may still be interesting to compute the latter instead of the former for computational efficiency; as we will illustrate in \S\ref{sec:timesNLP-IDvsNatural}, the natural bound can be computed faster than the $\NLPId$ bound.  

\begin{theorem} Consider the $\MESP$ instance $\MESP(C,s)$ and the related four problems \ref{fixinMi}, \ref{fixoutDi}, \ref{fixoutMj}, and\break  \ref{fixinDj}, with input matrices  $A$, $B$, $A_1$\,, $B_1$\,, $A_0$\,, $B_0$\,, defined as above. For all $i,j\in N$, we have 
\phantom{.}
\begin{enumerate}
    \item $\znatural(A_1,B_1,n-s) + (s-1)\log\lambda_{1}(C/C_{ii})  +  \log(C_{ii}) \leq \znatural(A_{(N \setminus i )\cdot},B,n-s) + s \log\lambda_{\max}\,.$
    \item $\znatural(A_0,B_0,n-1-s) + s\log\lambda_{1}(C_{N\setminus  j ,N\setminus  j }) \leq \znatural\left(A_{(N \setminus j )\cdot},\left({\begin{smallmatrix}\!\!A_{j\cdot}\!\!\\ B \end{smallmatrix}}\right),n-s-1\right) + s \log\lambda_{\max}\,.$
\end{enumerate}
\end{theorem}

\begin{proof}
\phantom{.}

\noindent 1.\quad Let $\hat{x}$ be a feasible  solution for $\NLPId$ related to the $\MESP$ instance $\MESP(\allowbreak C/C_{ii},s-1)$ with finite objective value. 
Then, 
    \begin{align*}
        \ldet&\left(A_{(N \setminus i )\cdot}^\top \Diag(\mathbf{e}-\hat x)A_{(N \setminus i )\cdot} + B^\top B\right)\\ &= \ldet\left(A_{(N \setminus i )\cdot}^\top \Diag(\mathbf{e}-\hat x)A_{(N \setminus i )\cdot} + A_{i\cdot}^\top A_{i\cdot} - A_{i\cdot}^\top A_{i\cdot} + B^\top B\right)\\
        &= \ldet\left(I_n -
\left({\begin{smallmatrix}\!\!A_{i\cdot}\!\!\\\!\! A_{(N \setminus i )\cdot}\!\!\end{smallmatrix}}\right)^\top
\Diag\left(\left(\begin{smallmatrix}
            1\\
            \hat x
\end{smallmatrix}\right)\right)
\left({\begin{smallmatrix}\!\!A_{i\cdot}\!\!\\\!\! A_{(N \setminus i )\cdot}\!\!\end{smallmatrix}}\right)
\right)\\
&= \ldet\left(I_n -\Diag\left(\left(\begin{smallmatrix}
            1\\
            \hat x
\end{smallmatrix}\right)\right)^{\scriptscriptstyle 1/2}
\left({\begin{smallmatrix}\!\!A_{i\cdot}\!\!\\\!\! A_{(N \setminus i )\cdot}\!\!\end{smallmatrix}}\right)\left({\begin{smallmatrix}\!\!A_{i\cdot}\!\!\\\!\! A_{(N \setminus i )\cdot}\!\!\end{smallmatrix}}\right)^\top
\Diag\left(\left(\begin{smallmatrix}
            1\\
            \hat x
\end{smallmatrix}\right)\right)^{\scriptscriptstyle 1/2}\right),
      \end{align*}
    where we use  that $A^\top A+B^\top B=I_n$\, in the second equation.
Let
\begin{align*}
    &H:= \begin{pmatrix}
    \underset{\scriptscriptstyle 1\times 1}{H_{11}} & \underset{\scriptscriptstyle 1\times (n-1)}{H_{12}}\\
    \underset{\scriptscriptstyle (n-1)\times 1}{H_{21}} & \underset{\scriptscriptstyle (n-1)\times (n-1)}{H_{22}}
\end{pmatrix}\\
&:= \begin{pmatrix}
    \frac{1}{\lambda_{\max}}C_{ii} & \frac{1}{\lambda_{\max}}C_{i,N\setminus  i }\Diag(\hat{x})^{\scriptscriptstyle 1/2}\\
    \frac{1}{\lambda_{\max}}\Diag(\hat{x})^{\scriptscriptstyle 1/2}C_{N\setminus  i ,i} & ~~~\Diag(\mathbf{e}-\hat x) +  \frac{1}{\lambda_{\max}} \Diag(\hat{x})^{\scriptscriptstyle 1/2} C_{N\setminus  i ,N\setminus  i }\Diag(\hat{x})^{\scriptscriptstyle 1/2}  
\end{pmatrix},
\end{align*}
and note that $H = I_n -\Diag\left(\left(\begin{smallmatrix}
            1\\
            \hat x
\end{smallmatrix}\right)\right)^{\scriptscriptstyle 1/2}
\left({\begin{smallmatrix}\!\!A_{i\cdot}\!\!\\\!\! A_{(N \setminus i )\cdot}\!\!\end{smallmatrix}}\right)\left({\begin{smallmatrix}\!\!A_{i\cdot}\!\!\\\!\! A_{(N \setminus i )\cdot}\!\!\end{smallmatrix}}\right)^\top
\Diag\left(\left(\begin{smallmatrix}
            1\\
            \hat x
\end{smallmatrix}\right)\right)^{\scriptscriptstyle 1/2}$. From the Schur determinant formula, we have
\begin{align*}
    &\ldet(H) = \ldet(H_{11}) + \ldet(H/H_{11})\\
    &= \log(\textstyle\frac{1}{\lambda_{\max}}C_{ii})  + \ldet\left(
    \vphantom{ +  (\textstyle\frac{1}{\lambda_{\max}}) \Diag(\hat{x})^{\scriptscriptstyle 1/2}\left(C_{N\setminus  i ,N\setminus  i } - C_{N\setminus  i ,i} C_{ii}^{-1}C_{i,N\setminus  i }\right)\Diag(\hat{x})^{\scriptscriptstyle 1/2}}
    \Diag(\mathbf{e}-\hat x)\right. \\
    &\qquad \left. +  (\textstyle\frac{1}{\lambda_{\max}}) \Diag(\hat{x})^{\scriptscriptstyle 1/2}\left(C_{N\setminus  i ,N\setminus  i } - C_{N\setminus  i ,i} C_{ii}^{-1}C_{i,N\setminus  i }\right)\Diag(\hat{x})^{\scriptscriptstyle 1/2}\right)\\
     &= \log(C_{ii}) - \log \lambda_{\max}  + \ldet\left(I_{n-1} +  \Diag(\hat{x})^{\scriptscriptstyle 1/2}(\textstyle\frac{1}{\lambda_{\max}} C/C_{ii} - I_{n-1} )\Diag(\hat{x})^{\scriptscriptstyle 1/2}\right).
\end{align*}
Then, we conclude that 
\begin{align*}
    &\ldet(A_{(N \setminus i )\cdot}^\top \Diag(\mathbf{e}-\hat x)A_{(N \setminus i )\cdot} + B^\top B) + s\log\lambda_{\max}\\
    &\quad= \ldet(I_{n-1} +  \Diag(\hat{x})^{\scriptscriptstyle 1/2}\left((\textstyle\frac{1}{\lambda_{\max}}) C/C_{ii} - I_{n-1} \right)\Diag(\hat{x})^{\scriptscriptstyle 1/2})\\
    &\qquad + \log(C_{ii}) + (s-1)\log\lambda_{\max}\,.
\end{align*}
Let $\tilde\lambda_{\max} := \lambda_{1}(C/C_{ii})$. From Lemma \ref{lem:nlp-natural-equiv-every-x}, we have that
\begin{align*}
        &\ldet(A_1^\top\Diag(\mathbf{e}-\hat{x})A_1 + B_1^\top B_1) + \log(C_{ii}) + (s-1)\log\tilde\lambda_{\max} \\
        &\quad= \ldet\left( I_{n-1} + \Diag(\hat{x})^{\scriptscriptstyle 1/2}\left(\textstyle\frac{1}{\tilde\lambda_{\max}^{}}C/C_{ii}-I_{n-1}\right)\Diag(\hat{x})^{\scriptscriptstyle 1/2}\right) \\
        &\qquad + \log(C_{ii})+(s-1)\log\tilde\lambda_{\max}\,.
    \end{align*}
     Noting that  $\lambda_{\max} \geq \tilde\lambda_{\max}$\, (see  \cite[Cor. 2.3]{SchurBook}), the result follows from Lemma \ref{lem:natural_to_nlp_modified}.

\medskip
    
\noindent 2. \quad Let $\hat{x}$ be a feasible  solution for $\NLPId$ related to the $\MESP$ instance $\MESP(\allowbreak C_{N\setminus  j ,N\setminus  j },s)$ with finite objective value. Let $\tilde{\lambda}_{\max} := \lambda_{1}(C_{N\setminus  j ,N\setminus  j })$.  From Lemma \ref{lem:nlp-natural-equiv-every-x}, we have
    \begin{align*}
        &\ldet\left(A_0^\top\Diag(\mathbf{e}-\hat{x})A_0 + B_0^\top B_0\right) + s\log\tilde\lambda_{\max} \\
        &\quad= \ldet\left( I_{n-1} + \Diag(\hat{x})^{\scriptscriptstyle 1/2}\left(\textstyle\frac{1}{\tilde\lambda_{\max}^{}} C_{N\setminus  j ,N\setminus  j }-I_{n-1}\right)\Diag(\hat{x})^{\scriptscriptstyle 1/2}\right)+s\log\tilde\lambda_{\max}\,,
    \end{align*}
    and 
    \begin{align*}
        &\ldet\left(A_{(N \setminus j )\cdot}^\top\Diag(\mathbf{e}-\hat{x})A_{(N \setminus j )\cdot} + A_{j\cdot}^\top A_{j\cdot} + B^\top B\right) + s\log\lambda_{\max}\\
        &\quad= \ldet\left( I_{n-1} + \Diag(\hat{x})^{\scriptscriptstyle 1/2}\left(\textstyle\frac{1}{\tilde\lambda_{\max}^{}} C_{N\setminus  j ,N\setminus  j }-I_{n-1}\right)\Diag(\hat{x})^{\scriptscriptstyle 1/2} \right) + s\log\lambda_{\max}\,,
    \end{align*}
    where we use the fact that 
    $A_{(N \setminus j )\cdot}A_{(N \setminus j )\cdot}^\top = I_{n-1} - (1/{\lambda_{\max}})C_{N\setminus  j ,N\setminus  j }$\,.
    Noting that  $\lambda_{\max} \geq \tilde\lambda_{\max}$\, (see  \cite[Thm. 4.3.8]{HJBook}), the result follows from Lemma \ref{lem:natural_to_nlp_modified}.
\qed \end{proof}


\section{Numerical experiments}\label{sec:numexp}~

\subsection{Experimental setup}
Some of the $\MESP$ relaxations addressed in \S\ref{sec:MESPbounds} rely on a scaling factor $\gamma$. For \ref{linx} and \ref{bqp_original}, we present the results with the optimal scaling factor; see \cite{Mixing}. For $\NLPDi$, the scaling factor $\gamma$ should be selected in the interval $[1/d_{\max}\,,1/d_{\min}]$. In our experiments, we tested 100 evenly spaced values for $\gamma$ in this interval and report the results for the best one.  

When comparing different upper bounds, we present the gaps defined by the difference between the upper bounds for $\MESP$ and 0/1 $\DOPT$, and the lower bounds computed by the local-search heuristics from  \cite{KLQ} and \cite{PonteFampaLeeMPB}, respectively. 

We used the solvers KNITRO (see \cite{KNITRO}) and MOSEK (see \cite{MOSEK}) to solve the convex optimization problems. 
More specifically, we used KNITRO for all relaxations except for \ref{bqp_original} and to solve the Trace subproblem in $\NLPTr$. We implemented our algorithms in Julia v1.11.3. 
We used the parameter settings for the solvers aiming at their best performance. Next, we summarize the settings that we employed, so that it is possible to reproduce our experiments.
We employed KNITRO 14.0.0 (via the Julia wrapper KNITRO.jl v0.14.4), using \texttt{CONVEX = true}, \texttt{FEASTOL} = $10^{-6}$ (feasibility tolerance), \texttt{OPTTOL} = $10^{-6}$ (absolute optimality tolerance),
\texttt{ALGORITHM} = 1 (Interior/Direct algorithm), 
\texttt{HESSOPT} = 6 (KNITRO computes a limited-memory quasi-Newton BFGS Hessian; we used the default value of 
\texttt{LMSIZE} = 10 limited-memory pairs stored when approximating the Hessian). 
We employed MOSEK 10.2.15 (via the Julia wrapper MOSEKTools.jl  v0.15.5).

We ran our experiments on ``zebratoo'', a
32-core machine (running Windows Server 2022 Standard):
two Intel Xeon Gold 6444Y processors running at 3.60GHz, with 16 cores each, and 128 GB of memory.


\subsection{Transferring bounds from MESP to pure 0/1 D-Opt}\label{subsec:pure}

We compare five bounds  for the pure 0/1 $\DOPT$ problem: linx,  factorization,  $\NLPId$,   $\NLPDi$, and  BQP, all computed for the associated $\MESP$ problem  produced by the map $\mathcal{M}$.  We recall that the  natural bound is equal to the bound obtained when we compute the $\NLPId$ bound for the associated MESP problem. 

This comparison  complements the results of \cite[Section 6.1]{PonteFampaLeeMPB}, where the natural bound for pure 0/1 $\DOPT$  instances was compared only with the linx  and the factorization bounds for the $\MESP$ instances produced by the map $\mathcal{M}$. We use the same approach to construct test-instances,  generating 
normally-distributed elements for the  $n\times m$ full column-rank  matrices $A$, with mean $0$ and standard deviation $1$. We set $n:=120$,  $m:= 40,60$, and we vary $s$.

\begin{figure}[!ht]%
    \centering
    \subfloat[$m=40$]{{\includegraphics[scale=0.30]{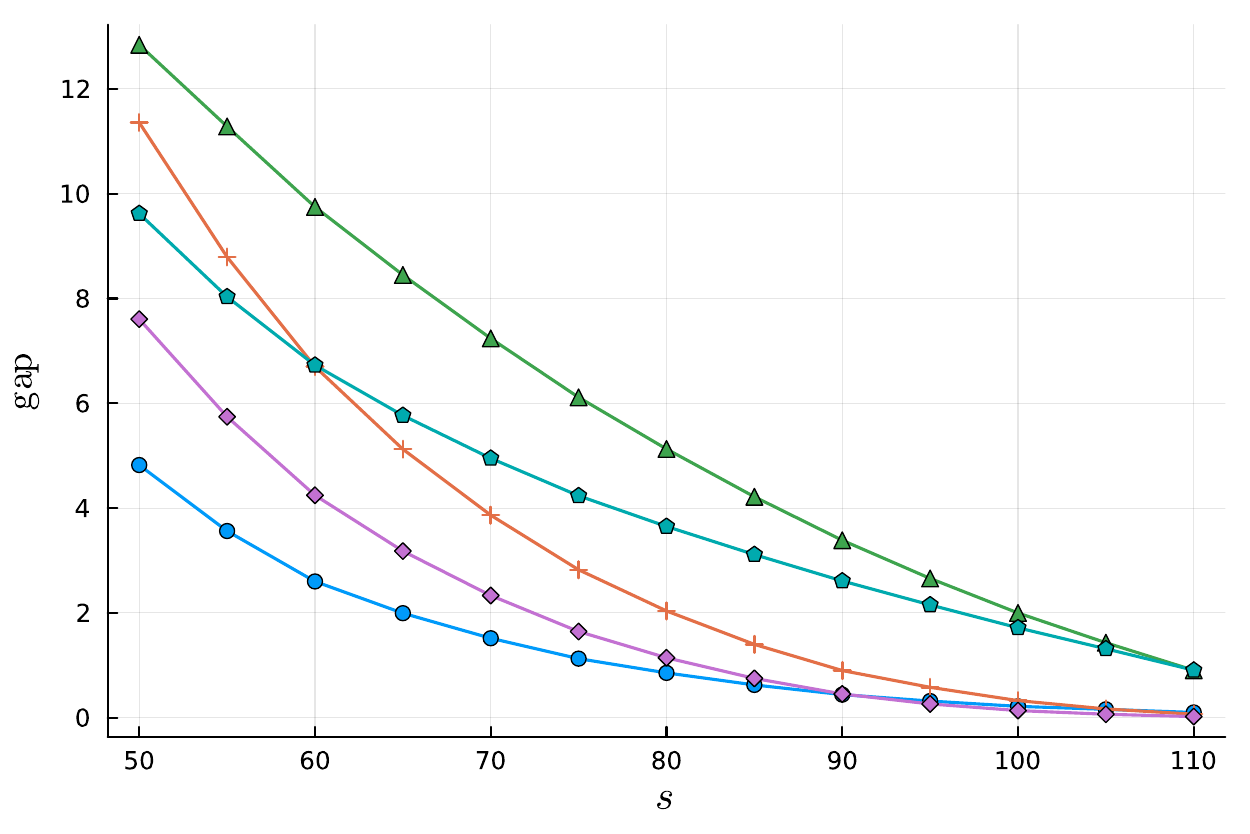} }}%
    \subfloat[$m=60$]{{\includegraphics[scale=0.30]{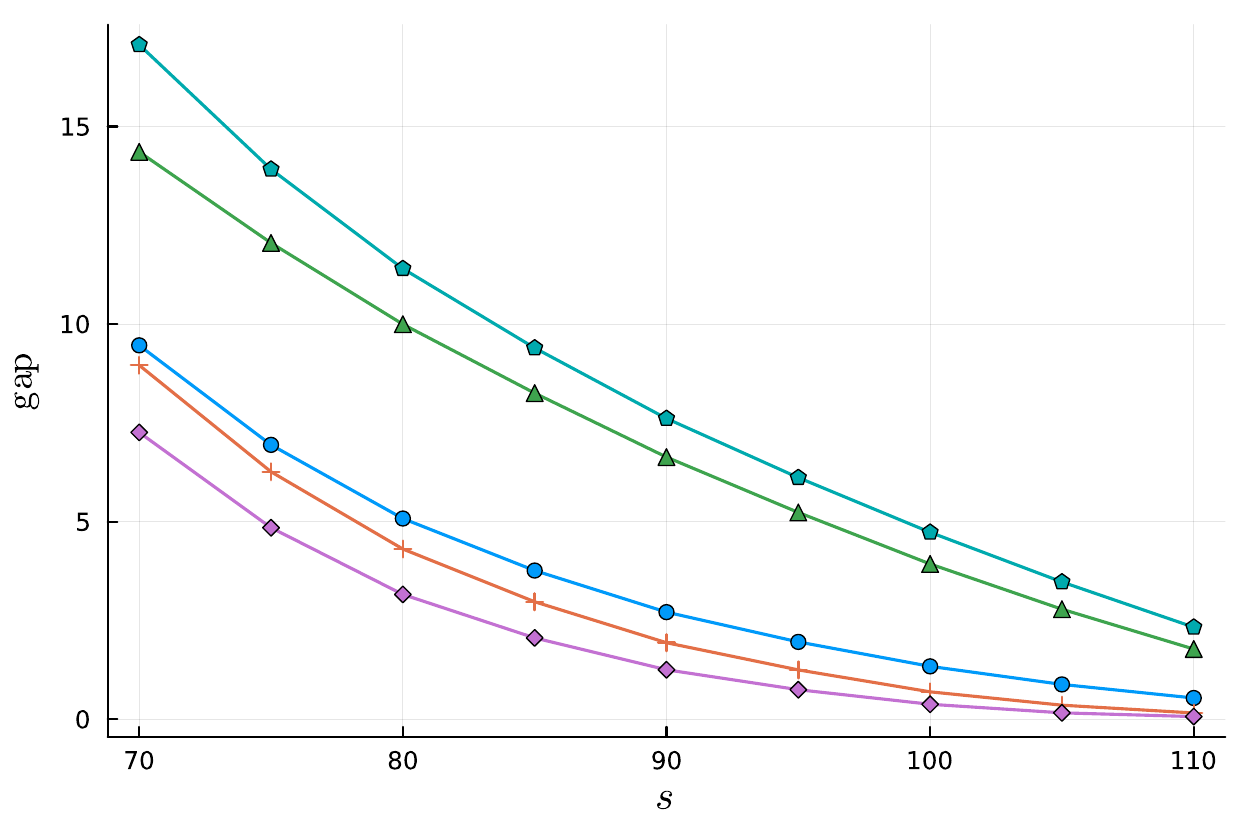} }}%
      \\
     \subfloat{{\includegraphics[scale=0.325]{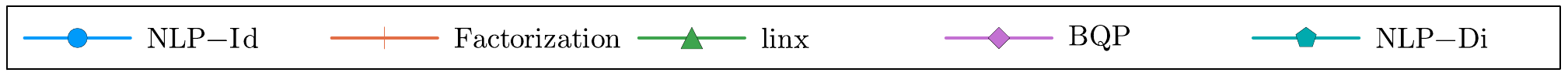} }}%
     \caption{$\mathcal{M}$-induced bounds for pure 0/1 $\DOPT$, varying $s$}%
    \label{fig:compare_bounds_randn}%
\end{figure}

In Figure \ref{fig:compare_bounds_randn}, we present plots comparing the bounds.  
The worst-performing bounds are linx and $\NLPDi$, which does \emph{not} follow the experience in the literature for positive-definite $\MESP$ instances; see,  e.g., \cite{AFLW_Using}, where $\NLPDi$ mainly outperforms  $\NLPId$, and \cite{FactPaper}, where linx is competitive with the factorization bound. We see that $\NLPId$ and \ref{bqp_original} are the best-performing bounds. Once more, the \ref{bqp_original} bound outperforming both the linx bound and the factorization bound was not expected after the results observed in \cite{gscale}. Interestingly, for these specific $\MESP$ instances,  the bounds that were not found to be  competitive with previous numerical results in the literature, namely, the $\NLPId$ bound and the BQP bound,  performed best.  The similar experiment in \cite[Sec. 6.1]{PonteFampaLeeMPB} also had poor performance for linx, but only compared to the factorization bound and the natural bound. 

\subsection{Impact of the multiplicity of the greatest eigenvalue of \texorpdfstring{$C$}{C} on MESP bounds}\label{subsec:multiplicity}

We can see comparing the two plots in Figure 
\ref{fig:compare_bounds_randn} that there appears to be a dependence in $m$ on both which bounds perform best and overall performance (smaller gaps with smaller $m$). We conducted an experiment to further investigate the behavior of bounds on $\MESP$ instances coming from 0/1 $\DOPT$,with varying $m$.
To construct the 0/1 $\DOPT$ instances, we randomly generated a $100\times 50$ full-column-rank matrix and then considered the $100\times m$  matrices $A$ as its submatrices, with the first $m:=50,45,\ldots,5$ columns. We set $s:=50$. 

\begin{figure}[!h]%
    \centering
    \subfloat[$\MESP$  from pure 0/1 $\DOPT$ instances] 
{{\includegraphics[width=.5\textwidth]{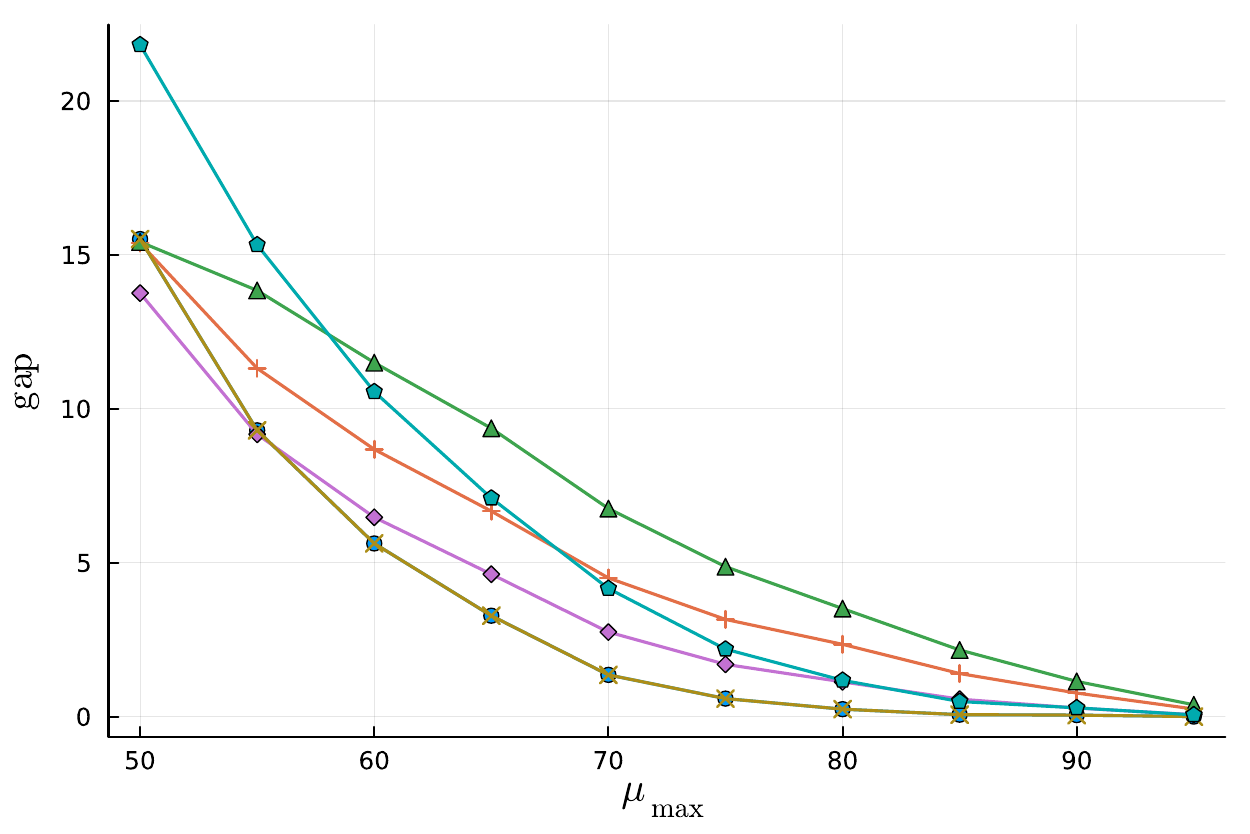} }}%
    \subfloat[\ref{MESP} with similar top eigenvalues]{{\includegraphics[width=.5\textwidth]{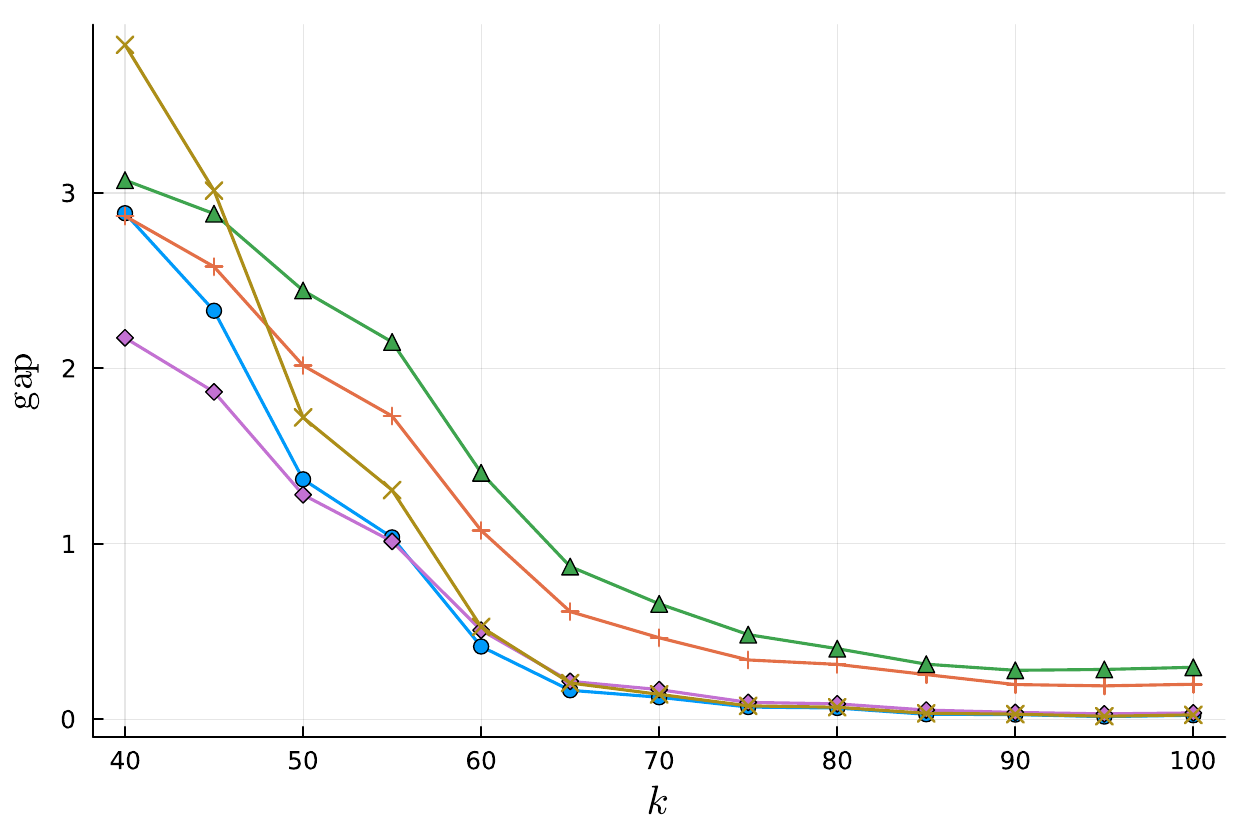} }}%
    \\
    \vspace{-5pt}\subfloat{{\includegraphics[width=.5\textwidth]{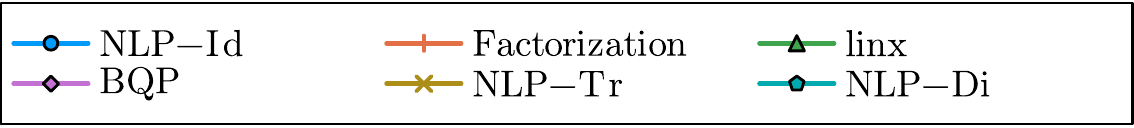} }}%
     \caption{$\MESP$ bounds for some particular $\MESP$ instances}%
     \label{fig:124}
\end{figure}

The results are in 
Figure \ref{fig:124}(a).
Note that the horizontal axis is labeled by 
$\mu_{\max}$\,, the multiplicity of the maximum eigenvalue of $A$.
We have that $\mu_{\max}=n-m$
($C=I_n-UU^\top$, where 
$U\Sigma V^\top$ is the compact singular-value decomposition of the 0/1 $\DOPT$ input matrix $A$. So $C$ has all $n-m$ positive eigenvalues equal to $1$).
Because our instances hold $n$ constant, decreasing $m$ is equivalent to increasing $\mu_{\max}$\,.
We compare the same $\MESP$ upper bounds as in Figure \ref{fig:compare_bounds_randn} and observe  similar results; for  $\mu_{\max}\leq 55$, the worst bound is $\NLPDi$, followed by linx, and the best bound is BQP, followed by $\NLPId$/$\NLPTr$. For $\mu_{\max}>55$, the worst bound is linx, and the best is $\NLPId$/$\NLPTr$. 
Finally, we see that all of the gaps in Figure \ref{fig:124}(a) are decreasing in $\mu_{\max}$\,, for which we have no convincing explanation. 

Next we sought to understand if our observations with regard to Figure \ref{fig:124}(a) are specific to 0/1 D-Opt. 
We consider  the singularity of $C$ and the multiplicity of its greatest eigenvalue.   
We modified a benchmark $n=124$ covariance matrix $C$.
We kept the eigenvectors of $C$, but set the $24$ least eigenvalues to zero, and modified the greatest $k-1$ eigenvalues of $C$
to be similar to the $k$-th greatest:
$\lambda_i(C) := 1.001^{k-i}\cdot \lambda_k(C)$, for $i = 1, \dots, k-1.$
 Our purpose was to investigate how similarity between the greatest eigenvalues (and singularity of $C$) would affect the behavior of the bounds. We set $s:=50$. We note for the (unaltered) benchmark instance with $n=124$ and $s=50$, the linx  and the factorization bounds are far superior to all of the others. 
The results are presented in Figure \ref{fig:124}(b). We do not show results for $\NLPDi$ because its performance was very poor. We confirm that the characteristics added to the benchmark instance alter the relationship between the bounds to the same behavior observed in Figure  \ref{fig:124}(a);   linx was the worst performing bound while BQP and $\NLPId$ performed best. 
Moreover, $\NLPId$ performs better than $\NLPTr$, which was also not observed for the benchmark instance. Finally, we see the fascinating phenomenon that all of the gaps in Figure \ref{fig:124}(b) are decreasing in $k$, following the behavior observed in Figure \ref{fig:124}(a).


\subsection{NLP-Id bound for \texorpdfstring{$\MESP(C,s)$}{MESP(C,s)} vs. the natural bound for
\texorpdfstring{$\mathcal{D}(\MESP(C,s);\Phi)$}{D(MESP(C,s);Phi)}}\label{sec:timesNLP-IDvsNatural}

We observed that when applying B\&B based on the $\NLPId$ bound to $\MESP$, we could, instead of computing the $\NLPId$ bound for the subproblems, compute the natural bound for the 0/1 $\DOPT$ instances derived from applying $\mathcal{D}$ to them. By Theorem \ref{thm:nlp-natural-equiv},  both approaches lead to  the same upper bound, but one of the bounds might be computed faster.  

For $C\in\mathbb{S}^n_+$ with real Schur decomposition $\Phi\Lambda\Phi^\top$, we consider $\DOPT(A, B,n-s):= \mathcal{D}(\mbox{MESP}(C,s);\Phi)$. We note that the $A$ matrix of the 0/1 $\DOPT$ instance is square. However, as we mentioned when we defined $\mathcal{D}$, 
$\rank(A)=n-\mu_{\max}$\,. So, it is possible to remove
$\mu_{\max}$ columns from 
{\scriptsize$\begin{pmatrix}A\\[-5pt] B \end{pmatrix}$}, 
making some upper-bound calculations for $\DOPT(A, B, n-s)$ much more efficient when $\mu_{\max}$ is large. Indeed, this is the case, as we establish with the following straightforward result.

\begin{proposition}\label{prop:compact_natural}
For  $C \in \mathbb{S}^n_{+}$\,, consider the $\MESP$ instance $\MESP(C,s)$, the real Schur decomposition $\Phi\Lambda\Phi^\top$ of $C$, and  $\mu_{\max}$\,,  the multiplicity of the greatest eigenvalue  of $C$.  Let $R := \{\mu_{\max}+1,\cdots,n\}$, and $\DOPT(A,B,n-s):= \mathcal{D}(\MESP(C,s);\Phi)$.  Then, we have  $\zrelax(A,B,n-s)=\zrelax(A_{\cdot R}\,,B_{\cdot R}\,,n-s)$, 
for $\mathcal{R}\in \{
\zdopt\,,\znatural\,,\zspectral\,,\zhadamard
\}$.
\end{proposition}

 From Proposition \ref{prop:compact_natural}, we see that the greater $\mu_{\max}$ is, the fewer the number of columns in the matrix $A$ in  the 0/1 $\DOPT$ instance derived from $\MESP(C,s)$ via the map $\mathcal{D}$ can be in our numerical calculations. We also observe that in the objective function of the natural bound, the argument of $\ldet(\cdot)$ is linear in $x$, which can be an additional advantage for computing the natural bound when compared to the $\NLPId$ bound.

\begin{figure}[!h]%
    \centering
    \subfloat[varying $s$, original matrix $C$]{{\includegraphics[width=.45\textwidth]{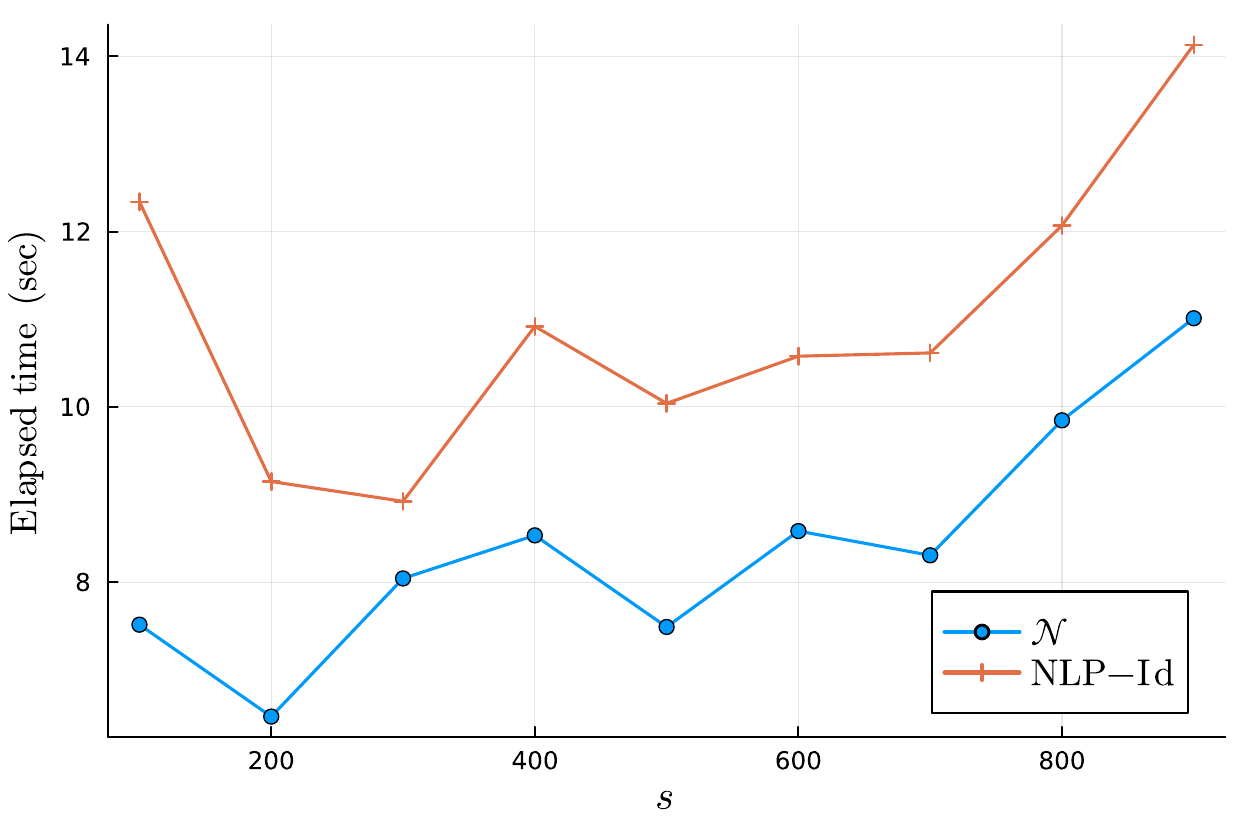} }}%
    \subfloat[varying $\mu_{\max}$ of $\tilde C$, $s:=900$]{{\includegraphics[width=.45\textwidth]{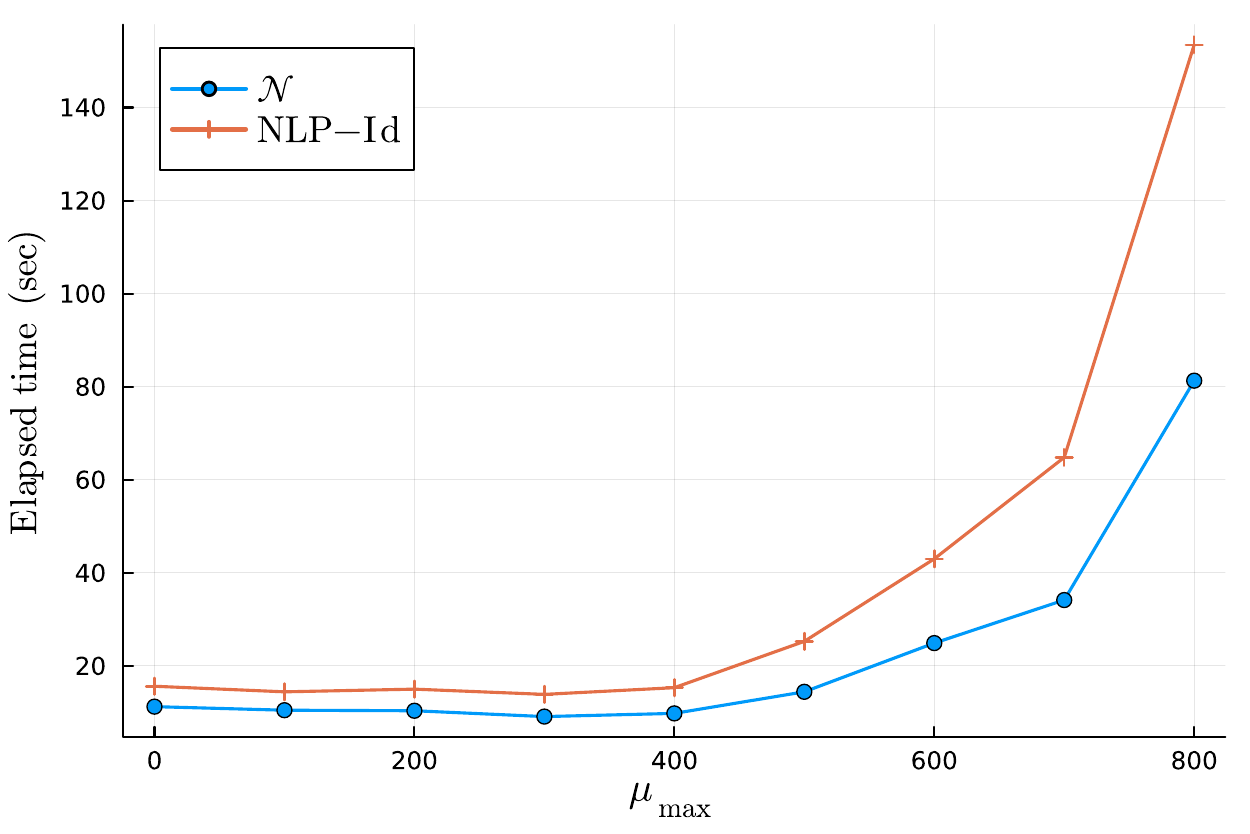} }}%
    \vspace{-7pt}
     \caption{Time comparison: NLP-Id bound vs. $\mathcal{D}$-induced natural bound}%
    \label{fig:compare_nlpid_natural}%
\end{figure}

Next, we compare the elapsed times for computing the natural bound and the $\NLPId$ bound for some  instances derived from a large-dimensional $\MESP$ instance from the literature; with the large dimension, we better illustrate the difference between the times. We consider an $n=2000$ covariance matrix  $C$ with rank 949 based 
on Reddit data from  \cite{Dey2018} and \cite{Munmun}. This instance was also used in experiments discussed in   
\cite{Weijun} and \cite{FactPaper}.  Initially, we consider the original covariance matrix $C$  with $\mu_{\max}=1$ and vary $s$ in $100,\dots,1500$. Next, we set $s:=900$ and construct new matrices $\tilde{C}$ with the same eigenvectors of $C$ and varying the multiplicity of its greatest eigenvalue by setting, for each $k\in \{1\}\cup \{100,200,\dots,800\}$,   $\lambda_i(\tilde{C})=\lambda_{k}(C)$, for $i=1,\ldots,k$, and $\lambda_i(\tilde{C})=\lambda_{i}(C)$, for $i=k+1,\ldots,n$. 
We depict the results in Figure \ref{fig:compare_nlpid_natural}. For all instances, we can observe a large improvement in time when solving the natural bound, with a significant benefit when $\mu_{\max}$ is large\,.

\FloatBarrier

\bibliographystyle{plain}
\bibliography{Bib}

\end{document}